\documentclass[12pt]{amsart}
\usepackage{amssymb, amsmath, amsthm}
\usepackage[utf8]{inputenc}
\usepackage{graphicx}
\usepackage{euscript}
\usepackage{dsfont}
\usepackage{color}
\usepackage{xcolor}
\usepackage{bm}
\usepackage{verbatim}

\definecolor{myred}{rgb}{0.2,0,0}
\definecolor{myblue}{rgb}{0,0,0.6}
\definecolor{mygreen}{rgb}{0,0.2,0}
\usepackage[hyphens]{url}
\usepackage[colorlinks=true,linkcolor=myred,citecolor=myblue,urlcolor=mygreen]{hyperref}
\usepackage[centering, bottom=1.5in]{geometry}
\usepackage[english]{babel}
\usepackage{bbm}
\usepackage{textcomp}
\usepackage{subfig}
\usepackage{float}

\usepackage[numbers,sort&compress]{natbib}

\newcommand{\R}{\mathbb{R}}

\newcommand{\Z}{\mathbb{Z}}

\newcommand{\norm}[1]{\left\lVert#1\right\rVert}

\renewcommand{\le}{\operatorname{\leqslant}}
\renewcommand{\ge}{\operatorname{\geqslant}}
\renewcommand{\leq}{\operatorname{\leqslant}}
\renewcommand{\geq}{\operatorname{\geqslant}}

\newtheorem*{conjecture*}{Conjecture}
\newtheorem{theorem}{Theorem}[section]
\newtheorem*{theorem*}{Theorem}
\newtheorem{lemma}{Lemma}[section]

\newtheorem{proposition}{Proposition}[section]

\newtheorem{remark}{Remark}[section]

\newtheorem{conjecture}{Conjecture}

\numberwithin{equation}{section}

\usepackage{caption}
\begin{document}
	\title{Poissonian pair correlation for $\alpha n^{\theta}$} 
	\author{Maksym Radziwi\l\l}
	\address{
          The University of Texas at Austin, Department of Mathematics, RLM 8.100, 2515 Speedway Stop C1200, Austin, TX 78712}
	\email{\href{mailto:maksym.radziwill@gmail.com}{maksym.radziwill@gmail.com}}
	
	\author{Andrei Shubin}
	\address{Institute of Discrete Mathematics and Geometry, TU Wien, Wiedner Hauptstr. 8-10, A-1040 Wien, Austria}
	\email{\href{mailto:andrei.shubin@tuwien.ac.at}{andrei.shubin@tuwien.ac.at}}
	
	\maketitle

	\begin{abstract}
	  We show that sequences of the form $\alpha n^{\theta} \pmod{1}$ with $\alpha > 0$ and $0 < \theta < \tfrac{43}{117} = \tfrac{1}{3} + 0.0341 \ldots$ have Poissonian pair correlation. This improves upon the previous result by Lutsko, Sourmelidis, and Technau, where this was established for $\alpha > 0$ and $0 < \theta < \tfrac{14}{41} = \tfrac{1}{3} +  0.0081 \ldots$.

          We reduce the problem of establishing Poissonian pair correlation to a counting problem using a form of amplification and the Bombieri-Iwaniec double large sieve. The counting problem is then resolved non-optimally by appealing to the bounds of Robert-Sargos and (Fouvry-Iwaniec-)Cao-Zhai. The exponent $\theta = \tfrac{2}{5}$ is the limit of our approach. 
                
        \end{abstract}

\section{Introduction}

A real-valued sequence $(x_n)_{n\ge 1}$ is \textit{equidistributed mod 1} if for any $0 < a < b < 1$, 
\[
	\lim_{N \to +\infty} \frac{1}{N} \# \bigl\{ n \le N:  \{ x_n \} \in [a, b]  \bigr\} = b-a.
\]
Equidistribution is one of the most basic properties of sequences that indicates their pseudorandom behavior.
Weyl~\cite{Weyl} showed that equidistribution modulo $1$ of a sequence $x_n$ is equivalent to showing cancellation in exponential sums
\[
	\sum_{n < N} e(k x_n) = o(N),
\]
for every fixed integer $k \neq 0$.
Subsequently, Weyl established that $\alpha n^d$ is equidistributed for all irrational $\alpha$ and $d \geq 1$ integer. Fejer and Csillag~(\cite[Corollary~2.1]{KN74}, \cite{Csillag}) established the same property when $d > 0$ is not an integer and $\alpha$ is non-zero.

Equidistribution is a ``large-scale'' property: it only provides information about the number of points inside intervals of size $\asymp 1$. A finer question is to ask about ``small-scale'' properties, that is distribution properties at the scale $\asymp \frac{1}{N}$. The most common statistic at this scale is the \textit{gap distribution} describing the spacings $N (y_{k + 1} - y_k)$ where $y_1 < y_2 < \ldots < y_N$ is an ordering of $x_1, \ldots, x_N \in [0,1]$.

It is natural to expect that the sequence $\alpha n^{\theta}$ with $\theta \not \in \mathbb{Z}$ and $\alpha \neq 0$ is distributed like a random set of points modulo $1$. This leads one to the following conjecture.
\begin{conjecture} \label{conj:main}
  Let $\alpha \neq 0$ and $\theta > 0$ with $\theta \neq \tfrac 12$ and $\theta \not \in \mathbb{Z}$. 
  Let $y_1^{(N)} < \ldots < y_N^{(N)}$ denote an ordering of the points $\alpha n^{\theta} \pmod{1}$ with $n < N$. Then,
  \[
  \lim_{N \rightarrow \infty} \frac{1}{N} \# \Big \{ i < N: N (y_{i + 1}^{(N)} - y_i^{(N)}) \in (a,b) \Big \} \rightarrow \int_{a}^{b} e^{-t} dt. 
  \]
  We then say that the gap distribution is Poissonian. 
\end{conjecture} 
We exclude the case $\theta = \tfrac 12$ because Elkies-McMullen~\cite{Elk_McM} have shown that the gap distribution \textit{is not} Poissonian in this case. This is also the only case in which the gap distribution is known to exist. A weaker statistic than the gap distribution is the pair correlation. The belief that the sequence $\alpha n^{\theta}$ behaves randomly leads to the following conjecture for the pair correlation. 
\begin{conjecture} \label{conj:pair}
  Let $\alpha \neq 0$ and $\theta > 0$ with $\theta \neq \tfrac 12$ and $\theta \not \in \mathbb{Z}$. Then, for any $s > 0$, 
  \[
  R_{2} \Big ( [-s, s], \alpha n^{\theta}, N \Big ) := \frac{1}{N} \cdot \# \Big \{ 1 \leq n \neq m \leq N: N \| \alpha n^{\theta} - \alpha m^{\theta} \| \leq s \Big \} \rightarrow 2 s
\]
as $N \rightarrow \infty$ and where $\| x \|$ denotes the distance of $x$ from the nearest integer. Whenever this holds we say that the sequence $\alpha n^{\theta}$ has Poissonian Pair Correlation (PPC). 
\end{conjecture}

One can think of the pair correlation as a second moment of the gap distribution. In particular if all higher $m$-correlations exist and are Poisson then Conjecture \ref{conj:main} would follow (see~\cite[Appendix~A]{Rudnick_Kurlberg}). 

We now briefly summarize the state of the art regarding Conjecture \ref{conj:pair}. There are broadly two types of results: metric results that allow for averaging in $\alpha$ or $\theta$, or deterministic results valid for specific $\alpha$ and $\theta$. As far as metric results go, PPC has been shown for sequences of the form $\alpha n^d$ with integer $d \ge 2$~\cite{Rudnick_Sarnak, Rudnick_Sarnak_Zaharescu, Mark_Strom, HB2010} for almost all $\alpha$. Moreover, Heath-Brown provided an algorithm for constructing a dense set of such $\alpha$'s in~\cite{HB2010}. The case $d=2$ is of particular interest due to its connection to quantum chaos (see~\cite{Berry_Tabor, Rudnick_08}). For sequences of the form $\alpha n^{\theta}$ with non-integer fixed $\theta > 0$ PPC is known for almost all $\alpha$ following Aistleitner, El-Baz, and Munsch~\cite{Aistl_Elb_Mun}, as well as Rudnick and Technau~\cite{Rudnick_Tech}. While in~\cite{Technau_Yesha} Technau and Yesha show that the sequence $n^{\theta}$ has PPC for almost all $\theta > 7$.

Deterministic results valid for specific $\theta$ and $\alpha$ are harder to come by. For $\theta = \tfrac 12$, El-Baz, Marklof, and Vinogradov~\cite{Elb_Mark_Vinog} showed that the sequence $(\sqrt{n})_{n \neq \square}$ admits PPC. Finally, Lutsko, Sourmelidis, and Technau~\cite{Lut_Sour_Tech} recently verified PPC for all sequences of the form $\alpha n^{\theta}$, where $\alpha > 0$ and $0 < \theta < \tfrac{14}{41} = \tfrac{1}{3} + 0.008\ldots$. Going beyond $\tfrac 13$ requires non-trivial bounds for certain exponential sums. The goal of this work is to extend the range of $\theta$ by deploying heavy exponential sums machinery.

\begin{theorem} \label{main_thm}
	The sequence $\alpha n^{\theta}$ has a Poissonian pair correlation for all $\alpha > 0$ and $0 < \theta < \tfrac{43}{117} = \tfrac{1}{3} + 0.03418\ldots$
\end{theorem} 

	\begin{figure}
			\centering
			\subfloat{{\includegraphics[width=8cm]{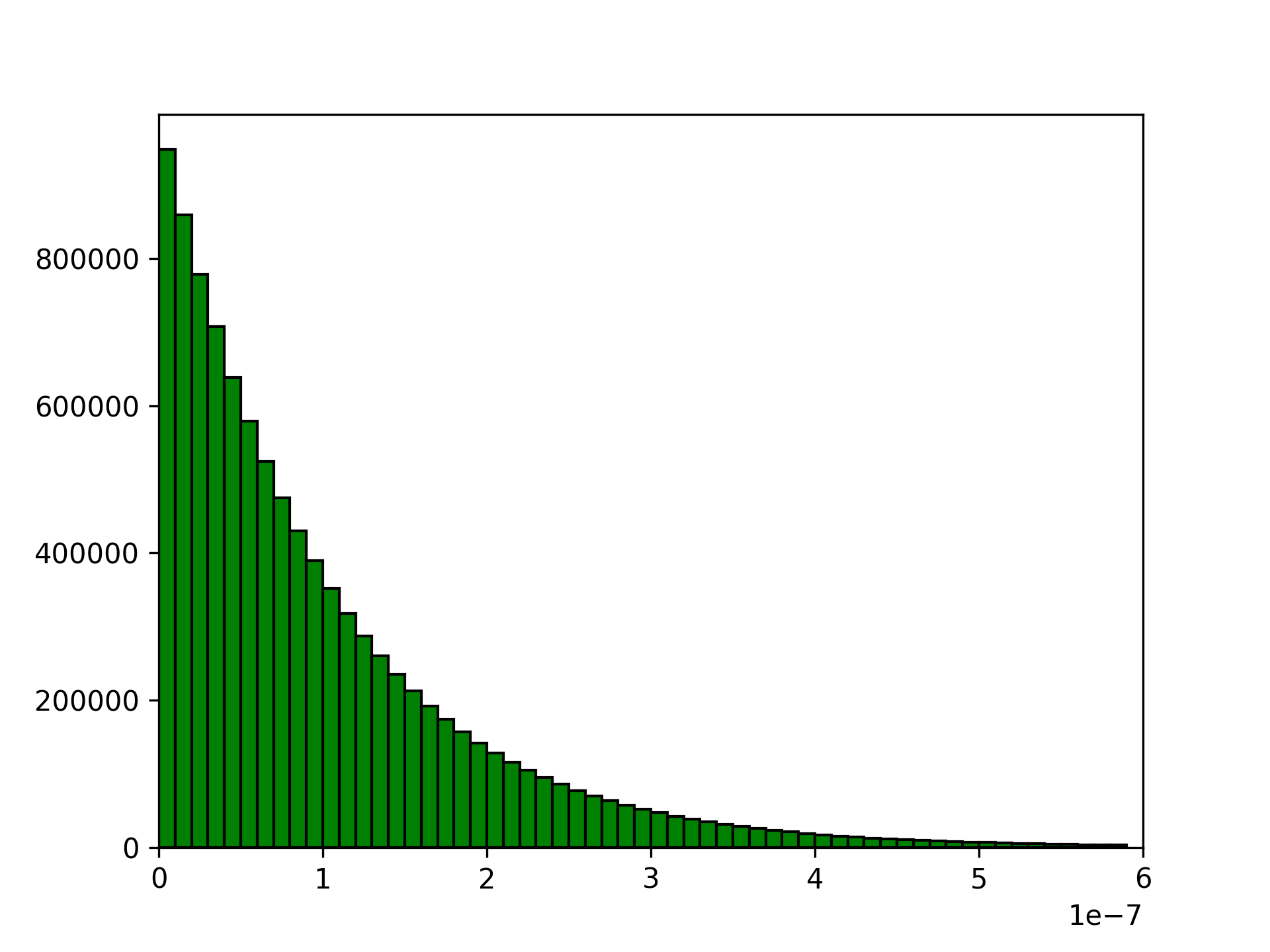}}}
			\subfloat{{\includegraphics[width=8cm]{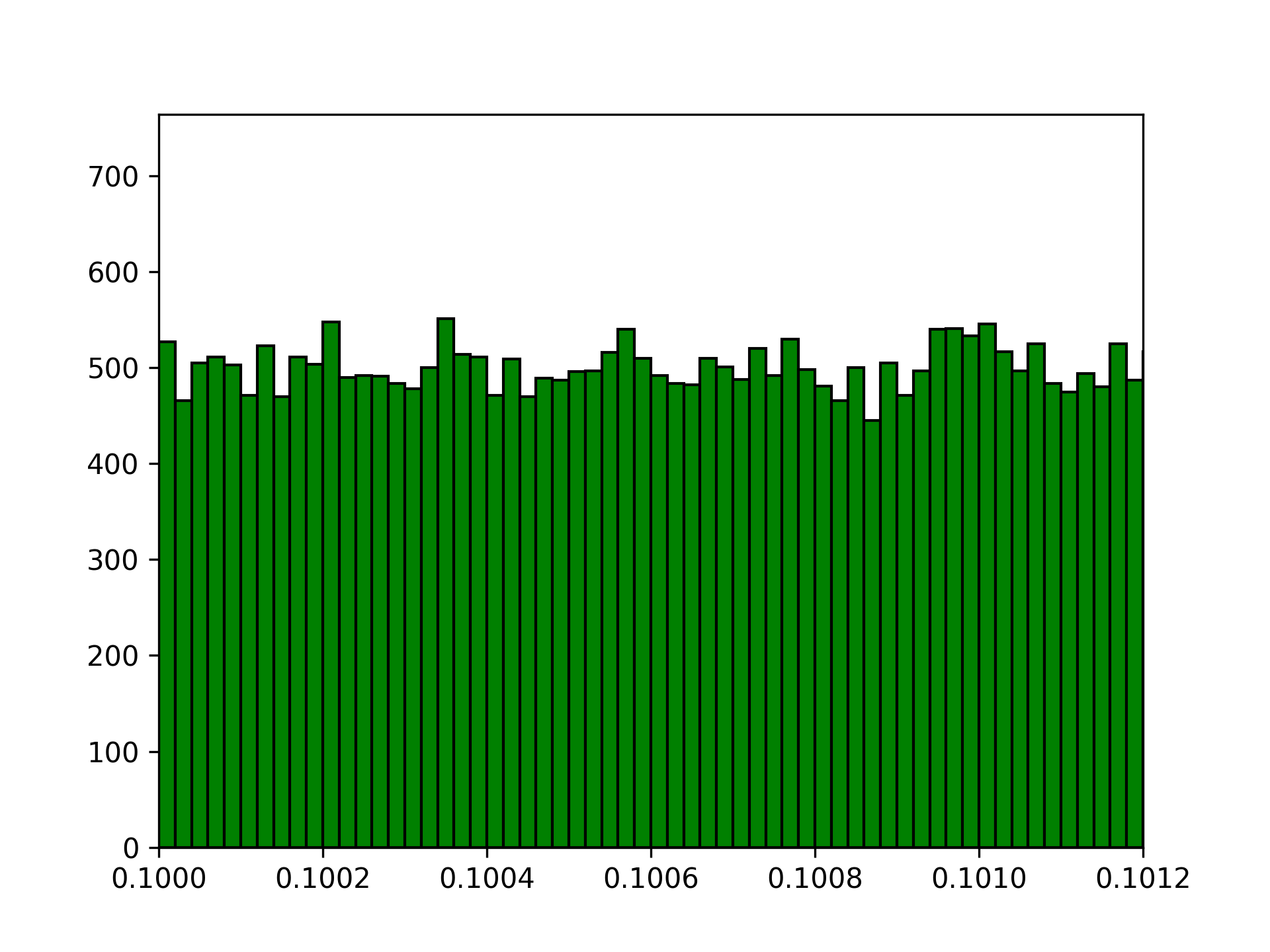} }}
			\caption{The distribution of gaps and spacings for the sequence $n^{1/3}$, with a step size $\frac{0.1}{N}$, $N = 10^7$ for gaps, and $N = 5000$ for spacings.}
	\end{figure}
	
	\begin{figure}
				\centering
				\subfloat{{\includegraphics[width=8cm]{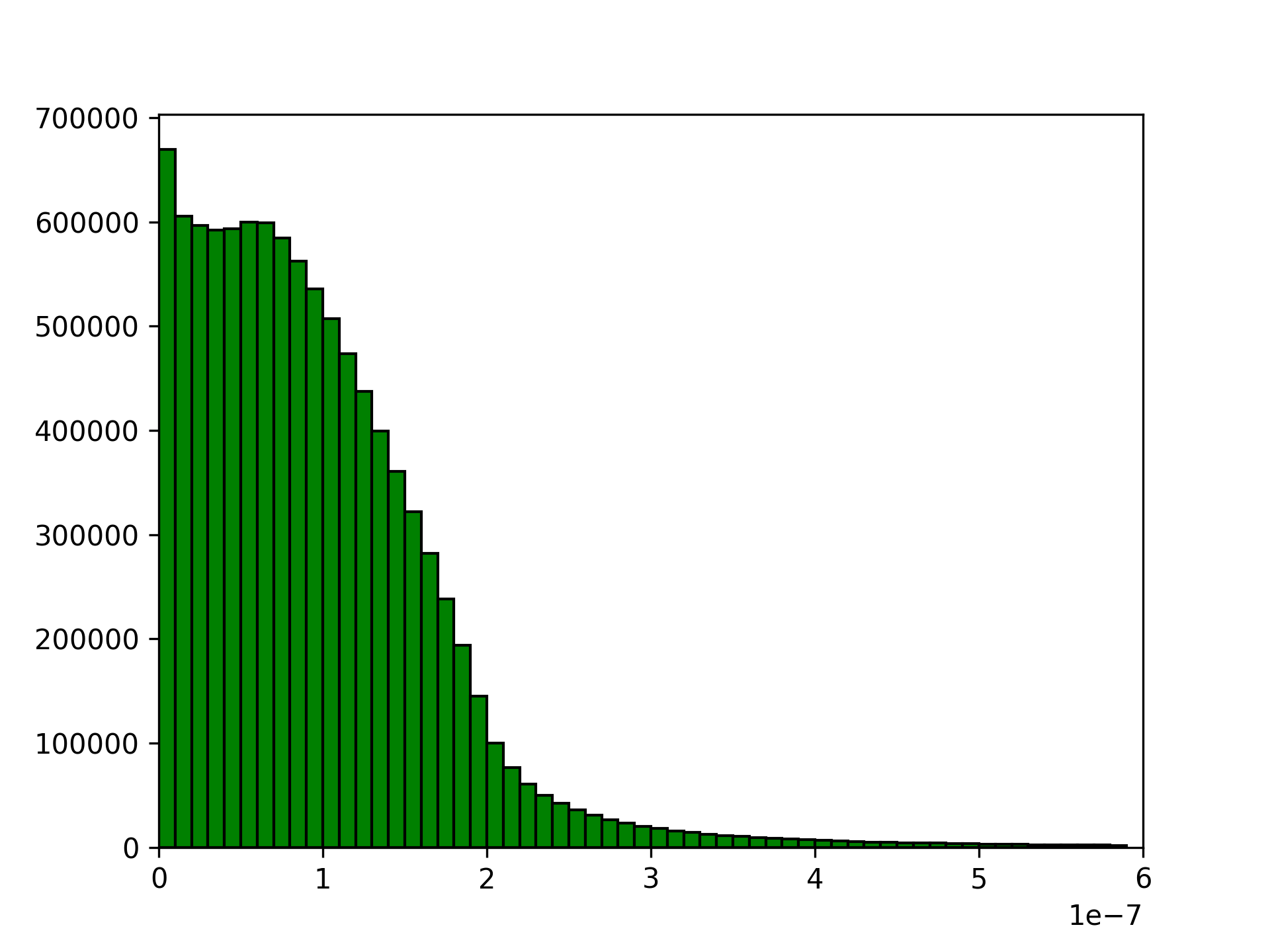}}}
				\subfloat{{\includegraphics[width=8cm]{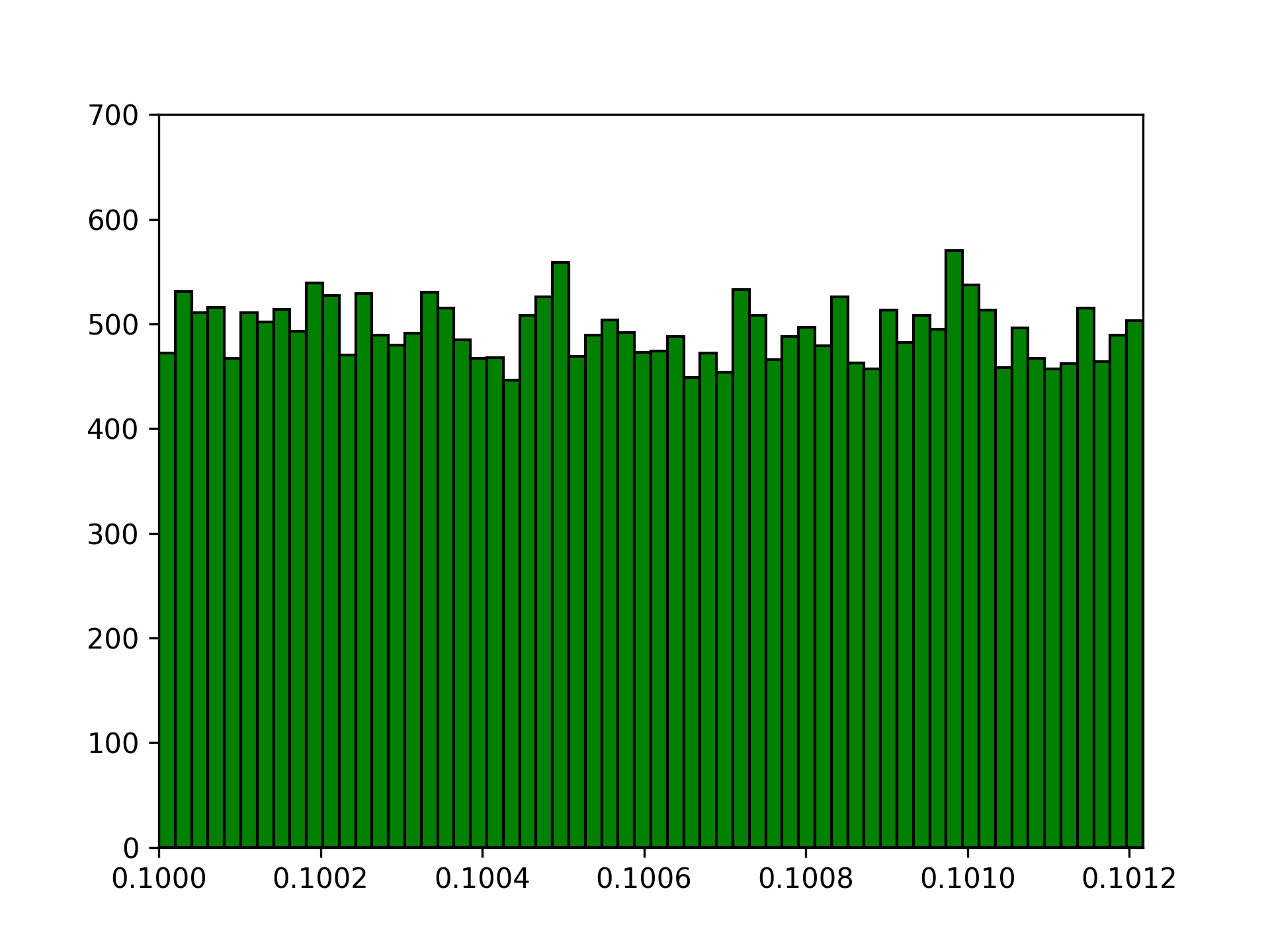} }}
				\caption{The distribution of gaps and spacings for the sequence $n^{1/2}$, with a step size $\frac{0.1}{N}$, $N = 10^7$ for gaps, and $N = 5000$, $n \neq \square$ for spacings. Elkies and McMullen~\cite{Elk_McM} demonstrated that the gap distribution in this case is not Poissonian, whereas El-Baz, Marklof, and Vinogradov~\cite{Elb_Mark_Vinog} showed that the pair correlation function is Poissonian.}
	\end{figure}
		
		\begin{figure}
					\centering
					\subfloat{{\includegraphics[width=8cm]{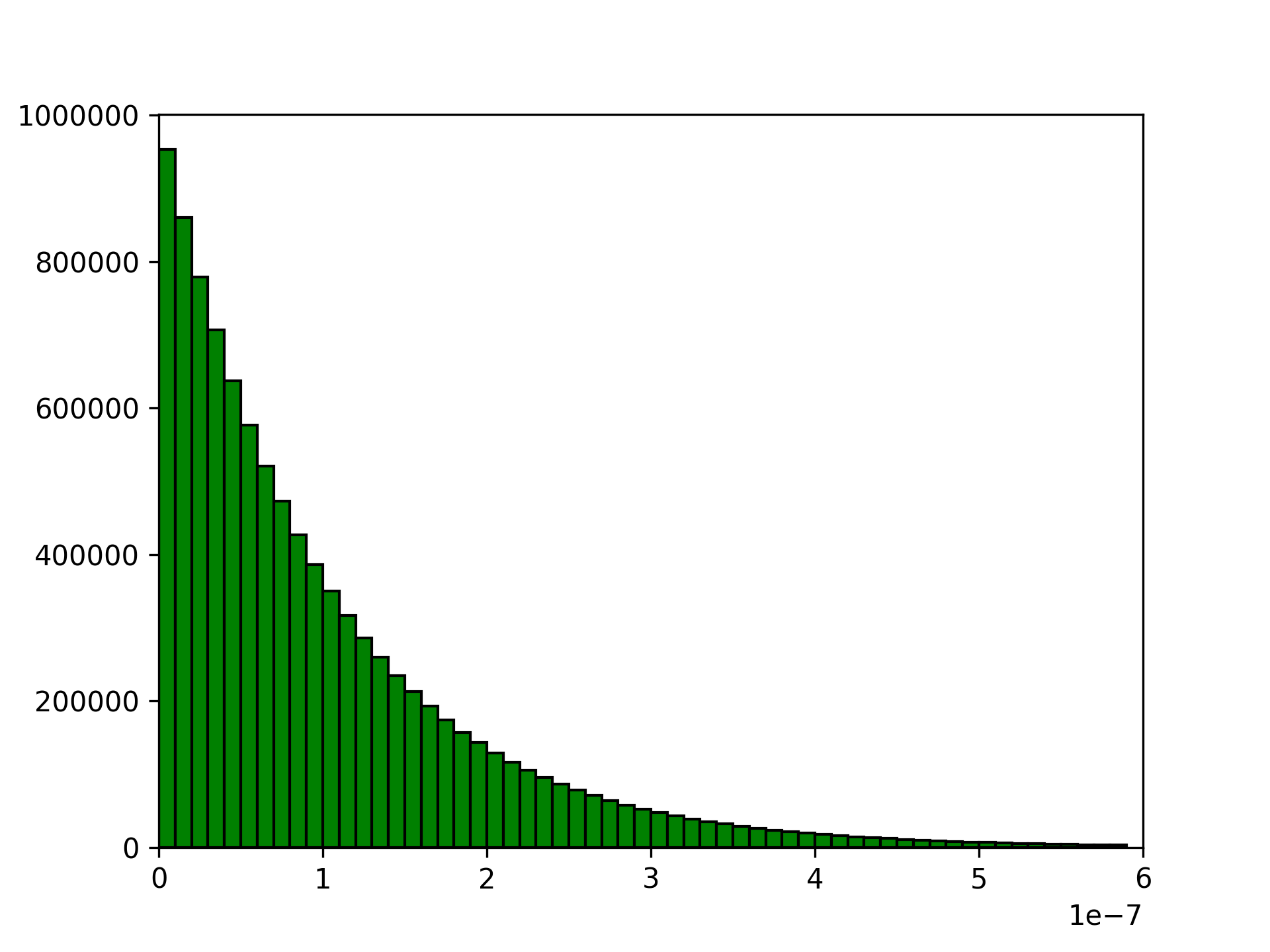}}}
					\subfloat{{\includegraphics[width=8cm]{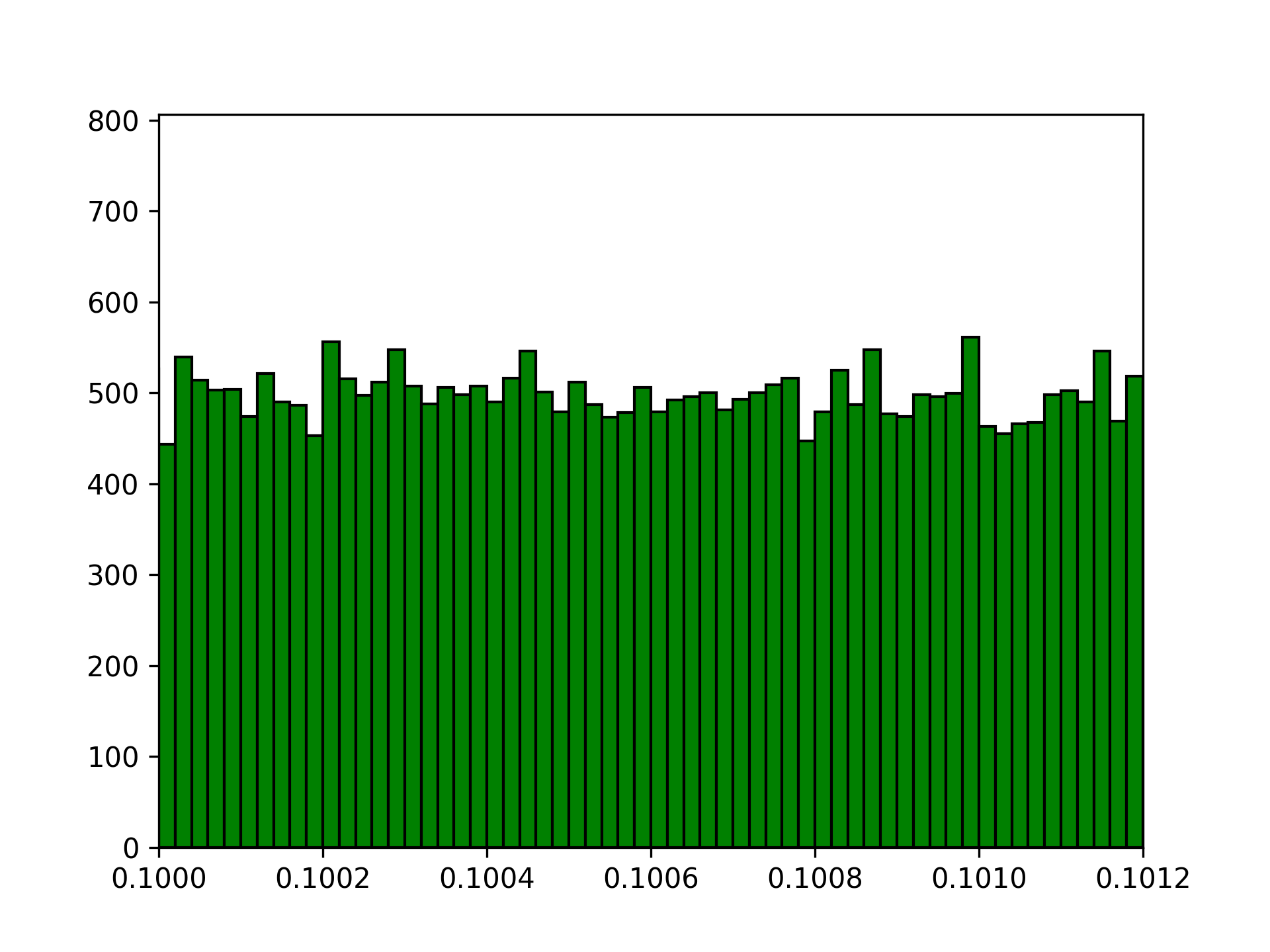} }}
					\caption{The distribution of gaps and spacings for the sequence $n^{2/3}$, with a step size $\frac{0.1}{N}$, $N = 10^7$ for gaps, and $N = 5000$ for spacings.}
		\end{figure}

\subsection{Reduction to exponential sums and sketch of the proof}

We define the integral version of the pair correlation function in the standard way:
\begin{equation} \label{eq:fourier}
	\tilde R_2 \bigl( f, \{ x_n \}, N \bigr) := \frac{1}{N} \sum_{k \in \Z} \sum_{1 \le i \neq j \le N} f\bigl( N(x_i - x_j + k) \bigr).
\end{equation} The condition $\norm{x_i - x_j} \le \tfrac{s}{N}$ is equivalent to $N(x_i - x_j + k) \in [-s, s]$ for some $k \in \Z$. When $N$ becomes large compared to $s$, there is only one such $k$. Then
\[
	\lim_{N \to + \infty} R_2 \bigl( [-s, s], \{ x_n \}, N \bigr) = \lim_{N \to +\infty} \tilde R_2 \bigl( \mathbbm{1}_{[-s, s]}, \{ x_n \}, N \bigr).
        \]
        Smoothing $\mathbbm{1}_{[-s,s]}$ and expanding into Fourier series reduces the problem of PPC for $\alpha n^{\theta}$ to a problem about exponential sums. We state this below.

        \begin{proposition} \label{prop:reduction}
          Let $\alpha > 0$ and $\theta > 0$ be given.  
          Suppose that for every even smooth compactly supported function $f$, every $\varepsilon > 0$,
          \[
       \lim_{N \rightarrow \infty} \frac{2}{N^2} \sum_{0 < |k| \leq N^{1 + \varepsilon}} \widehat{f} \Big ( \frac{k}{N} \Big ) \Big | \sum_{1 \leq y \leq N} e(\alpha k y^{\theta}) \Big |^2 = f(0).
       \]
       Then the sequence $\alpha n^{\theta}$ has Poissonian Pair Correlation (PPC). 
        \end{proposition}

        Thus Conjecture \ref{conj:pair} amounts to showing that certain exponential sums have square-root cancellation on average. 

Since the growth rate of the phase $\alpha k y^{\theta}$ is $\ll \alpha N^{1 + \theta + \varepsilon}$, its first derivative is much smaller than $N$. Therefore, a natural first step is to apply the Poisson summation formula to each of the three sums over $k$, $y_1$, and $y_2$ as it will shorten each of these sums to length $N^{\theta}$. This step is formalized in the following proposition.

\begin{proposition} \label{prop_1}
	Let $\alpha, \theta, \varepsilon$ be real numbers such that $\alpha > 0$, $\varepsilon > 0$, $3\varepsilon < \theta < 1 - 3\varepsilon$, and for $K, Y_1, Y_2 > 0$ define the subset $\mathcal{N} (K, Y_1, Y_2)$ as follows:
	\begin{gather*}
		\mathcal{N} (K, Y_1, Y_2) := \bigl\{ \eta: \eta = m_1^{-\theta / (1-\theta)} - m_2^{-\theta / (1-\theta)}, \text{ where} \\
		\alpha \theta K (2Y_1)^{\theta-1} \le m_1 < 2 \alpha \theta K Y_1^{\theta-1}; \\  \alpha \theta K (2Y_2)^{\theta-1} \le m_2 < 2 \alpha \theta K Y_2^{\theta-1}; \\
		m_1 < m_2 \bigr\}.
	\end{gather*} Then we have
	\begin{multline} \label{main_transform}
		\frac{2}{N^2} \sum_{k \le N^{1+\varepsilon}} \hat{f} \bigl( \frac{k}{N} \bigr) \Big| \sum_{1 \le y \le N} e\bigl( \alpha k y^{\theta} \bigr) \Big|^2 = \\
		f(0) + O\Big( N^{-\theta/2 + \varepsilon/2} \log N  + N^{-1/4 + \theta/4 + 3\varepsilon/4} \log N + S (N) \Big),
	\end{multline} where
	\begin{equation} \label{new_expsum}
		S(N) := \frac{1}{N^2} \sum_{\substack{K \le N^{1+\varepsilon} \\ \text{dyadic}}} \sum_{\substack{Y_1, Y_2 \le N \\ \text{dyadic}}} K^{(2-\theta) / (1-\theta) - 1} \max_{K \le \tilde K \le 2K} \bigl| \tilde S(K, \tilde K, Y_1, Y_2) \bigr|,
	\end{equation} 
	\begin{multline*}
		\tilde S(K, \tilde K, Y_1, Y_2) := \\
		\sum_{\eta \in \mathcal{N}(K, Y_1, Y_2)} \Big( \frac{1}{m_1 m_2} \Big)^{(2-\theta) / (2-2\theta)}
		\sum_{L \le \ell < \tilde L} \frac{c_1}{\sqrt{\eta}} \Big(\frac{\ell}{\eta}\Big)^{1/(2\theta) - 1} e\bigl( c_2\ell^{1/\theta} \eta^{1 - 1/\theta} \bigr),
	\end{multline*}
	\begin{gather*}
		L = \frac{c_3 \eta}{1-\theta} K_1^{\theta / (1-\theta)}, \qquad \tilde L = \frac{c_3 \eta}{1-\theta} K_2^{\theta / (1-\theta)}, \qquad c_3 = (\alpha \theta)^{1/(1-\theta)} \Big( \frac{1}{\theta} - 1 \Big), \\
		c_1 = \frac{1 - \theta}{\sqrt{c_3 \theta}} \Big( \frac{1-\theta}{c_3} \Big)^{1/(2\theta) - 1}, \qquad c_2 = -\theta \Big( \frac{1-\theta}{c_3} \Big)^{1/\theta - 1},
	\end{gather*}
	\begin{gather}
		K_1 := \max \Big( K, \frac{m_1}{2\alpha \theta Y_1^{\theta-1}},  \frac{m_2}{2\alpha \theta Y_2^{\theta-1}} \Big), \label{K_1} \\
		K_2 := \min \Big( \tilde K, \frac{m_1}{\alpha \theta (2Y_1)^{\theta-1}}, \frac{m_2}{\alpha \theta (2Y_2)^{\theta-1}} \Big). \label{K_2}
	\end{gather} The implied constant in the error term in~\eqref{main_transform} depends on $f$, $\alpha$, $\theta$ and $\varepsilon$. 
\end{proposition}

A similar transform was used in~\cite{Lut_Sour_Tech}. An application of Proposition~\ref{prop_1}, combined with a trivial upper bound for $|\tilde S(K, \tilde K, Y_1, Y_2)|$ on the right side of~\eqref{new_expsum} yields Theorem~\ref{main_thm} in the range $\varepsilon < \theta < \tfrac 13 - \varepsilon$. Thus, going beyond $\tfrac 13$ requires a non-trivial bound for $S(N)$. We note incidentally that such non-trivial bounds are not always possible, for example, we have $S(N) \gg 1$ for $\alpha = 1$ and $\theta \in \{\tfrac 12, 1, 2, 3, \ldots\}$. 

To a first approximation,
\[
S(N) \approx N^{-1/2 - 3\theta / 2} \sum_{\eta \in \mathcal{N}} \sum_{\ell \sim N^{\theta}} e \Big (\ell^{1 / \theta} \eta^{1 - 1 / \theta} \Big ) 
\]
where
\[
\mathcal{N} := \Big \{ m_1, m_2 \sim N^{\theta} : m_1^{- \theta / (1 - \theta)} - m_2^{- \theta / (1 - \theta)} \Big \}. 
\]
For the purpose of the sketch we can also assume that elements of $\mathcal{N}$ are in ``generic position'' and thus of size $\asymp N^{-\theta^2 / (1 - \theta)}$. Given $\eta \in \mathcal{N}$ we let $\widetilde{\eta} := \eta N^{\theta^2 / (1 - \theta)}$ so that for generic $\eta \in \mathcal{N}$ we have $|\widetilde{\eta}| \asymp 1$. Likewise for $\ell \sim N^{\theta}$ we denote by $\widetilde{\ell} := \ell / N^{\theta}$ so that $|\widetilde{\ell}| \asymp 1$.  
To bound $S(N)$ we apply Holder's inequality.
\[
S(N) \ll N^{-1/2 - 3 \theta/2} \Big ( \# \mathcal{N} \Big )^{1 - 1/k} \cdot \Big ( \sum_{\eta \in \mathcal{N}} \Big | \sum_{\ell \sim N^{\theta}} e \Big ( \ell^{1/\theta} \eta^{1 - 1 / \theta} \Big ) \Big |^{k} \Big )^{1/k}.
\]
We notice that the phase is of size $X := N^{1 + \theta}$. Using the double large sieve we have the bounds,
\begin{multline*}
 \sum_{\eta \in \mathcal{N}} \Big | \sum_{\ell \sim N^{\theta}} e \Big ( \ell^{1 / \theta} \eta^{1 - 1 / \theta} \Big ) \Big |^{k} \\ \ll X^{1/2} \cdot \Big ( \sum_{\substack{\eta_1, \eta_2 \in \mathcal{N}}} \mathbf{1} \Big ( \Big | {\widetilde{\eta}_1}^{1 - 1 / \theta} - {\widetilde{\eta}_2}^{1 - 1 / \theta} \Big | \leq \frac{1}{X} \Big ) \Big )^{1/2} \cdot \\
 \Big ( \sum_{\substack{\ell_1, \ldots, \ell_{2k} \sim N^{\theta}}} \mathbf{1} \Big ( \Big |\widetilde{\ell}_1 + \ldots - \widetilde{\ell}_{2k} \Big | \leq \frac{1}{X} \Big ) \Big )^{1/2}. 
\end{multline*}
In the best case scenario the probability of each event inside the indicator function is $1/X$. Taking care to exclude the diagonal, this gives in the best case scenario the bound,
\[
\sum_{\eta \in \mathcal{N}} \Big | \sum_{\ell \sim N^{\theta}} e \Big ( \ell^{1 / \theta} \eta^{1 - 1 / \theta} \Big ) \Big |^{k} \ll N^{(1 + \theta)/2} \Big ( N^{2\theta} + \frac{N^{4 \theta}}{N^{1 + \theta}} \Big )^{1/2} \cdot \Big ( N^{k \theta} + \frac{N^{2k \theta}}{N^{1 + \theta}} \Big )^{1/2}.
\]
If $\theta_{k}$ is the optimal exponent that such an approach can yield, then we find $\theta_{2} = \tfrac 13$, $\theta_{3} = \theta_4 = \tfrac{2}{5}$, $\theta_{5} = \tfrac{5}{13}$, $\theta_{6} = \tfrac{3}{8}$ followed by smaller exponents. Thus the most beneficial choice of exponent is $k = 4$. Unfortunately we are not able to directly resolve the problem of counting eight-tuples of $(\ell_1, \ldots, \ell_8)$ all of size $\approx N^{\theta}$ such that,
\[
|\ell_1^{1 / \theta} + \ldots - \ell_{8}^{1 / \theta} | \leq N^{-\theta}. 
\]
Instead we perturb the problem a little and apply van der Corput differencing twice, to get that,
\begin{multline*}
	\sum_{\eta \in \mathcal{N}} \Big | \sum_{\ell \sim N^{\theta}} e (\ell^{1 / \theta} \eta^{1 - 1 / \theta}) \Big |^4 \ll \\
	\frac{N^{4 \theta} \# \mathcal{N}}{H_1^2} + \frac{N^{4 \theta} \# \mathcal{N}}{H_2} + \frac{N^{3 \theta}}{H_1 H_2} \sum_{\eta \in \mathcal{N}} \sum_{\substack{0 < |h_1| \leq H_1 \\ 0 < |h_2| \leq H_2}} \gamma(h_1, h_2) \sum_{\ell \sim N^{\theta}} e \Big (\eta^{1 - 1 / \theta} \ t (\ell, h_1, h_2) \Big )
\end{multline*}
where
\[
t (\ell, h_1, h_2) := ((\ell + h_1 + h_2)^{1 / \theta} - (\ell + h_2)^{1 / \theta}) - ((\ell + h_1)^{1 / \theta} - \ell^{1 / \theta}).
\]
We now apply the double large sieve, and this leads to the bound
\[
\sum_{\eta \in \mathcal{N}} \sum_{\substack{0 < |h_1| \leq H_1 \\ 0 < |h_2| \leq H_2}} \gamma(h_1, h_2) \sum_{\ell \sim N^{\theta}} e \Big (\eta^{1 - 1 / \theta} t (\ell, h_1, h_2) \Big ) \ll X^{1/2} \mathcal{B}_1^{1/2} \mathcal{B}_2^{1/2}
\]
where $X := N^{1 - \theta} H_1 H_2$ and
\[
\mathcal{B}_1 := \sum_{\substack{\eta_1, \eta_2 \in \mathcal{N}}} \mathbf{1} \Big ( | {\widetilde{\eta}_1}^{1 - 1 / \theta} - {\widetilde{\eta}_2}^{1 - 1 / \theta} | \leq \frac{1}{X} \Big )
  \]
  and
  \[
\mathcal{B}_2 := \sum_{\substack{0 < |h_1|, |h_1'| \leq H_1 \\ 0 < |h_2|, |h_2'| \leq H_2}} \sum_{\substack{\ell, \ell' \sim N^{\theta}}} \mathbf{1} \Big ( \Big | \widetilde{t}(\ell, h_1, h_2) - \widetilde{t}(\ell', h_1', h_2') \Big | \leq \frac{1}{X} \Big ) 
  \]
where as before $\widetilde{t}(\ell, h_1, h_2)$ denotes a suitably normalized version of $t(\ell, h_1, h_2)$ so that $|\widetilde{t}(\ell, h_1, h_2)| \asymp 1$. We bound $\mathcal{B}_1$ using a Taylor expansion and the bound of Robert-Sargos \cite{Robert-Sargos}, which gives an optimal bound for the number of tuples $(m_1, m_2, m_3, m_4)$ of size $M$ such that
\[
|m_1^{\alpha} + m_2^{\alpha} - m_3^{\alpha} - m_4^{\alpha} | \leq \delta.
\]
We then bound $\mathcal{B}_2$ using results of Cao and Zhai~\cite{Cao-Zhai}. Such bounds for ``shifted variables'' go back to the work of Fouvry-Iwaniec \cite{Fourvy_Iwaniec} with subsequent improvements by Liu~\cite{Liu} and Sargos-Wu~\cite{Sargos_Wu}. Due to the complicated nature of the bound of Cao-Zhai (see Lemma \ref{Cao-Zhai}) we refrain from stating the exact bounds here. It suffices to say that after inputting all these bounds we find an appropriate choice of $H_1$ and $H_2$ that yields the following main technical proposition:

\begin{proposition} \label{prop_2}
	Let $\alpha, \varepsilon, \mathcal{N}$ and $S(N)$ be as in Proposition~\ref{prop_1}, and let $\theta$ be such that $3\varepsilon < \theta < \tfrac{43}{117} - 5\varepsilon$. Then we have
	\[
		S(N) = o(1) \qquad \text{when } N \to +\infty,
	\] where the implied constant depends on $\alpha$, $\theta$ and $\varepsilon$. 
\end{proposition}

Theorem~\ref{main_thm} clearly follows from Propositions~\ref{prop:reduction},~\ref{prop_1} and~\ref{prop_2}. We prove Proposition~\ref{prop:reduction} in Section~\ref{sec_prop_reduction}. We then prove Proposition~\ref{prop_1} in Section~\ref{sec_pf_prop_1}, and Proposition~\ref{prop_2} in Section~\ref{sec_pf_prop_2}.

\subsection{Notation}

We use the standard notation for the real character $e(x) := e^{2 \pi i x}$. 

The relations $f(x) \ll g(x)$ and $f(x) = O(g(x))$ mean $|f(x)| \le C g(x)$ for a fixed number $C > 0$ and all large enough $x$. The symbols $\ll_{\varepsilon}$ and $O_{\varepsilon}$ mean that the constant $C$ might depend on the parameter $\varepsilon$. The relation $f(x) = o(g(x))$ for $g(x) > 0$ means $f(x) / g(x) \to 0$ as $x \to +\infty$. The relation $f(x) \asymp g(x)$ means $f(x) \ll g(x) \ll f(x)$. Finally, the relation $f(x) \sim g(x)$ means $f(x) / g(x) \to 1$ as $x \to +\infty$. 

For brevity, the limits in the sums over dyadic intervals, such as $\sum_{K, Y_1, Y_2, R}$, are often not indicated.

\subsection{Acknowledgments}

The first author received support from NSF grant DMS-1902063. The second author acknowledges the support provided by the Austrian Science Fund FWF (I4945-N). The second author also thanks Christoph Aistleitner, Tim Browning, and members of their research groups for the helpful conversations during his short visits.

\section{Reduction to exponential sums}
\label{sec_prop_reduction}
We prove in this section Proposition \ref{prop:reduction}.
   There is a sequence of $C_c^{\infty} (\R)$ functions $f_j^{+}$ and $f_j^{-}$, such that $f_j^{-} \le \mathbbm{1}_{[-s, s]} \le f_j^{+}$ for all $j$,
\[
	\int_{\R} \bigl( f_j^{+}(x) - f_j^{-} (x) \bigr) dx \to 0
\] when $j \to +\infty$, and $\hat{f_j^{\pm}} (x) \ll x^{-A}$ for arbitrarily large $A > 0$. Then, for the proof of Theorem~\ref{main_thm}, it is enough to show that for any $ f \in C_c^{\infty} (R)$, one has
\[
	\tilde R_2 \bigl( f, \{ \alpha n^{\theta} \}, N \bigr) \to \int_{\R} f(x) dx
\] when $N \to +\infty$. Applying Poisson summation to the sum over $k \in \Z$, we obtain
\[
	\tilde R_2 \bigl(f, \{ \alpha n^{\theta} \}, N \bigr) = \frac{1}{N^2} \sum_{k \in \Z} \hat{f} \Big( \frac{k}{N} \Big) \sum_{1 \le y_1 \neq y_2 \le N} e\bigl( \alpha k (y_1^{\theta} - y_2^{\theta} ) \bigr).
\] Moving away the term $k = 0$, we find
\[
	\tilde R_2 \bigl(f, \{ \alpha n^{\theta} \}, N \bigr) = \frac{1}{N^2} \sum_{\substack{k \in \Z \\ k \neq 0}} \hat{f} \Big( \frac{k}{N} \Big) \sum_{1 \le y_1 \neq y_2 \le N} e\bigl( \alpha k (y_1^{\theta} - y_2^{\theta} ) \bigr) + \hat{f}(0) + O\Big( \frac{\hat{f}(0)}{N} \Big).
\] For convenience we add and subtract the diagonal term $y_1 = y_2$: 
\begin{multline*}
	\tilde R_2 \bigl(f, \{ \alpha n^{\theta} \}, N \bigr) = \\
	\frac{1}{N^2} \sum_{\substack{k \in \Z \\ k \neq 0}} \hat{f} \Big( \frac{k}{N} \Big) \Big| \sum_{1 \le y \le N} e\bigl( \alpha k y^{\theta} \bigr) \Big|^2 - \frac{1}{N^2} \sum_{k \in \Z} \sum_{1 \le y \le N} \hat{f} \Big( \frac{k}{N} \Big) + \hat{f}(0) + O\Big( \frac{1}{N} \Big).
\end{multline*} Note that the term $k=0$ is absorbed by the error $O(1/N)$. By Poisson summation,
\[
	\frac{1}{N^2} \sum_{k \in \Z} \sum_{1 \le y \le N} \hat{f} \Big( \frac{k}{N} \Big) = \frac{1}{N} \sum_{k \in \Z} \hat{f} \Big( \frac{k}{N} \Big) = f(0).
\] Thus,
\[
	\tilde R_2 \bigl(f, \{ \alpha n^{\theta} \}, N \bigr) = \frac{1}{N^2} \sum_{\substack{k \in \Z \\ k \neq 0}} \hat{f} \Big( \frac{k}{N} \Big) \Big| \sum_{1 \le y \le N} e\bigl( \alpha k y^{\theta} \bigr) \Big|^2 - f(0) + \hat{f}(0) + O\Big( \frac{1}{N} \Big).
\] Since $\hat{f}(0) = \int_{\R} f(x) dx$, the problem reduces to showing that
\[
	\frac{1}{N^2} \sum_{\substack{k \in \Z \\ k \neq 0}} \hat{f} \Big( \frac{k}{N} \Big) \Big| \sum_{1 \le y \le N} e\bigl( \alpha k y^{\theta} \bigr) \Big|^2 = f(0) + o(1). 
\] Without loss of generality, we can assume that $f$ is even. Then the last formula is equivalent to
\[
	\frac{2}{N^2} \sum_{k=1}^{+\infty} \hat{f} \Big( \frac{k}{N} \Big) \Big| \sum_{1 \le y \le N} e\bigl( \alpha k y^{\theta} \bigr) \Big|^2 = f(0) + o(1). 
\] Due to the fast decay of $\hat{f} (x)$, the summation over $k$ can be restricted to $k \le N^{1+\varepsilon}$ for arbitrary small $\varepsilon > 0$. Thus for the proof of Theorem~\ref{main_thm} it is enough to show that
\begin{equation} \label{original_sum}
	\frac{2}{N^2} \sum_{k \le N^{1+\varepsilon}} \hat{f} \Big( \frac{k}{N} \Big) \Big| \sum_{1 \le y \le N} e\bigl( \alpha k y^{\theta} \bigr) \Big|^2 = f(0) + o_{\varepsilon} (1)
\end{equation} for any $\theta < \tfrac{43}{117} - 5\varepsilon$ (say).

\section{Proof of Proposition~\ref{prop_1}}
\label{sec_pf_prop_1}

\subsection{Poisson summation on $y$}

We use the following form of the Poisson summation formula (see~\cite[Theorem~8.16]{Iwaniec-Kowalski}):

\begin{lemma} \label{Poisson_summation}
	Let $f(x)$ be a real function on $[a,b]$ such that
	\[
		\Lambda \le f''(x) \le c \Lambda, \qquad |f^{(3)}(x)| \le c \Lambda (b-a)^{-1}, \qquad |f^{(4)}(x)| \le c \Lambda (b-a)^{-2}
	\] for some $\Lambda > 0$ and $c \ge 1$. Then
	\[
		\sum_{a < n < b} e\bigl( f(n) \bigr) = \sum_{\alpha < m < \beta}  \frac{1}{\sqrt{f''(x_m)}} e\Big( f(x_m) - mx_m + \frac{1}{8} \Big) + R_f (a,b), 
	\] where $\alpha = f'(a), \beta = f'(b)$, $f'(x_m) = m$,
	\[
		R_f (a,b) \ll \Lambda^{-1/2} + c^2 \log (\beta - \alpha + 1).
	\]
\end{lemma}

When $Y_1$ and $Y_2$ are much different in size we estimate their contribution using van der Corput's theorem (see~\cite[Theorem~8.20]{Iwaniec-Kowalski}):
\begin{lemma} \label{Corput_kth_test}
	Let $b-a \ge 1$, $f(x)$ be a real function on $[a,b]$, and $\nu \ge 2$ such that $\Lambda \le f^{(\nu)}(x) \le c\Lambda$, where $\Lambda > 0$ and $c \ge 1$. Then
	\[
		\sum_{a < n < b} e\bigl( f(n) \bigr) \ll \Lambda^{\kappa} (b-a) + \Lambda^{-\kappa} (b-a)^{1 - 2^{2-\nu}}  
	\] where $\kappa = (2^{\nu} - 2)^{-1}$ and the implied constant is absolute. 
\end{lemma}

Let us split the sum over $y \le N$ in~\eqref{original_sum} to dyadic intervals $(Y, 2Y]$, so that $Y \le N/2$. Then, applying Lemma~\ref{Poisson_summation} with $f(y) = -\alpha k y^{\theta}$, substituting $m \to -m$, and taking the conjugation, we get
\[
	\sum_{Y \le y < 2Y} e\bigl( \alpha k y^{\theta} \bigr) = T_k (Y) + R_k (Y),
\] where
\begin{gather} \nonumber
	T_k (Y) = \sum_{\alpha \theta k (2Y)^{\theta - 1} \le m < \alpha \theta k Y^{\theta-1}} \frac{c_4}{\sqrt{k}} \Big( \frac{k}{m} \Big)^{(2-\theta) / (2-2\theta)} e\Big( c_3 k^{1/(1-\theta)} m^{-\theta/(1-\theta)}  - \frac{1}{8} \Big), \\ \nonumber
	c_4 = \frac{1}{\sqrt{\alpha \theta (1-\theta)}} \bigl( \alpha \theta \bigr)^{(2-\theta) / (2-2\theta)}, \qquad c_3 = \bigl( \alpha \theta \bigr)^{1/(1-\theta)} \Big( \frac{1}{\theta} - 1 \Big), \\ \label{R_k_bound}
	R_k (Y) \ll \frac{Y^{1 - \theta/2}}{\sqrt{k}} + \max \bigl( 1, \log(k Y^{\theta-1}) \bigr). 
\end{gather} Similarly, split the sum over $k \le N^{1+\varepsilon}$ in~\eqref{original_sum} to dyadic intervals $(K, 2K]$:
\begin{multline*}
	\frac{2}{N^2} \sum_K \sum_{K \le k < 2K} \hat{f} \bigl( \frac{k}{N} \bigr) \Big| \sum_Y \bigl( T_k (Y) + R_k (Y) \bigr) \Big|^2 =: \\ 
	\frac{2}{N^2} \sum_K \Big( S_{1,1} (K) + S_{1,2} (K) + S_{2,1} (K) + S_{2,2} (K) \Big),
\end{multline*} where the sums $S_{1,1}, S_{1,2}, S_{2,1}, S_{2,2}$ correspond respectively to $T_k (Y_1) \overline{T_k (Y_2)}$, $T_k (Y_1) \overline{R_k (Y_2)}$, $R_k (Y_1) \overline{T_k (Y_2)}$, $R_k (Y_1) \overline{R_k (Y_2)}$. \\

We start by estimating the contribution from $S_{2,2}$. Applying trivially $\hat{f}~(k/N)~\ll~1$ and~\eqref{R_k_bound}, we find
\[
	S_{2,2}(K) = \sum_{K \le k < 2K} \hat{f} \Big( \frac{k}{N} \Big) \sum_{Y_1, Y_2} R_k (Y_1) \overline{R_k (Y_2)} \ll \\ K (\log N)^2 \Big( \Big( \frac{N^{1 - \theta/2}}{\sqrt{K}} \Big)^2 + \bigl( \log K \bigr)^2 \Big).    
\] Hence,
\begin{equation} \label{S_2_2_final_bound}
	\frac{2}{N^2} \sum_K S_{2,2}(K) \ll N^{-\theta} (\log N)^3 + N^{-1+\varepsilon} (\log N)^5,
\end{equation} which is negligible. \\

The sums $S_{1,2}$ and $S_{2,1}$ are similar, so we only consider $S_{1,2}$:
\[
	S_{1,2}(K) = \sum_{K \le k < 2K} \hat{f} \Big( \frac{k}{N} \Big) \sum_{Y_1, Y_2} T_k (Y_1) \overline{R_k (Y_2)}.
\] Applying Cauchy inequality to the sum over $K \le k < 2K$, we obtain
\begin{multline*}
	S_{1,2} (K) \le \sum_{Y_1, Y_2} \Big( \sum_{K \le k < 2K} \bigl| \hat{f} \Big( \frac{k}{N} \Big) \bigr|^2 |R_k (Y_2)|^2 \frac{k^{(2-\theta) / (1-\theta)}}{k} \Big)^{1/2} \cdot \\
	\Big( \sum_{K \le k < 2K} \Big| \sum_{\alpha \theta k (2Y_1)^{\theta - 1} \le m < \alpha \theta k Y_1^{\theta-1}} \Big( \frac{c_4}{m} \Big)^{(2-\theta) / (2-2\theta)} e\Big( c_3 k^{1 / (1-\theta)} m^{-\theta / (1-\theta)} - \frac{1}{8} \Big) \Big|^2 \Big)^{1/2}.
\end{multline*} Hence,
\begin{equation} \label{S_1_2}
	S_{1,2} (K) \ll \sum_{Y_1, Y_2} \Big( \sum_{K \le k < 2K} |R_k (Y_2)|^2 k^{1/(1-\theta)} \Big)^{1/2} \bigl| U_{1,2} (K)\bigr|^{1/2},
\end{equation} where
\[
	U_{1,2} (K) :=
	\sum_{K \le k < 2K} \sum_{\substack{m_1, m_2 \in \\ [\alpha \theta k (2Y_1)^{\theta-1}, \alpha \theta k Y_1^{\theta-1})}} \Big( \frac{c_4^2}{m_1 m_2} \Big)^{(2-\theta) / (2-2\theta)} e\Big(  c_3 k^{1/(1-\theta)} \eta(m_1, m_2) \Big),
\] where $\eta(m_1, m_2) = m_1^{-\theta / (1-\theta)} - m_2^{-\theta / (1-\theta)}$. Changing the order of summation we obtain
\[
	U_{1,2} (K) = \sum_{\substack{m_1, m_2 \in \\ [\alpha \theta K (2Y_1)^{\theta-1}, 2\alpha \theta K Y_1^{\theta-1})}} \Big( \frac{c_4^2}{m_1 m_2} \Big)^{(2-\theta) / (2-2\theta)}
	\mathop{{\sum}^{'}}_{K \le k < 2K} e\Big(  c_3 k^{1/(1-\theta)} \eta(m_1, m_2) \Big),
\] where $\mathop{{\sum}^{'}}$ denotes the sum with the additional restrictions 
\[
	m_1, m_2 \in [ \alpha \theta k (2Y_1)^{\theta-1}, \alpha \theta k Y_1^{\theta-1} ].
\]

Next, substitute
\[
	m_1 := m, \qquad m_2 := m+r.
\] Without loss of generality, let us assume that $m_1 \le m_2$, which implies that $r \ge 0$. Next, we can split the sum $U_{1, 2}$ into its diagonal ($r=0$) and non-diagonal ($r > 0$) terms. Furthermore, we can split the non-diagonal sum into dyadic intervals $r \in (R, 2R]$, where $R$ satisfies $R \le (2 - 2^{\theta-1}) \alpha \theta K Y_1^{\theta - 1}$. We get
\begin{multline} \label{U_1_2_split}
	U_{1,2} (K) = D_{1,2} (K) + \\
	\sum_R \sum_{R \le r < 2R} \sum_{\substack{m: \ m, m+r \in \\ [\alpha \theta K (2Y_1)^{\theta-1}, 2\alpha \theta K Y_1^{\theta-1})}} \Big( \frac{c_4^2}{m(m+r)} \Big)^{(2-\theta) / (2-2\theta)}  \mathop{{\sum}^{'}}_{K \le k < 2K} e\Big( c_3 k^{1/(1-\theta)} \eta \Big),
\end{multline} where $\eta = \eta(m,m+r) = m^{-\theta / (1-\theta)} - (m+r)^{-\theta / (1-\theta)} \sim r m^{-1/(1-\theta)}$. Note that the sum over $k$ can be empty, in which case the necessary bound follows immediately. 

\begin{remark}
	From now on, we assume that $\alpha \theta K Y_1^{\theta-1} \ge 10$. In this case, the sum over $m$ contains more than one term. Otherwise, we trivially have $U_{1,2} (K) \ll K$, which, by~\eqref{S_1_2}, implies
	\[
		S_{1,2} (K) \ll \sum_{\substack{Y_1, Y_2 \\ Y_1 \gg K^{1/(1-\theta)}}} \Big( \Big( \frac{Y_2^{1-\theta/2}}{\sqrt{K}} \Big)^2 K^{1/(1-\theta)} \Big)^{1/2} K^{1/2} \ll K^{1/(2-2\theta)} N^{1-\theta/2} \log N,
	\] which contributes to the right hand side of~\eqref{main_transform} no more than
	\[
		\frac{2}{N^2} \sum_{K \ll N^{1-\theta}} K^{1/(2-2\theta)} N^{1-\theta/2} \log N \ll N^{-1/2 - \theta/2} \log N, 
	\] which is negligible.
\end{remark}

For the diagonal term $D_{1,2} (K)$, we have
\begin{equation} \label{diagonal}
	D_{1,2} (K) \le \sum_{K \le k < 2K}\sum_{\substack{m \in \\ [\alpha \theta K (2Y_1)^{\theta-1}, 2\alpha \theta K Y_1^{\theta-1})}} \Big( \frac{c_4}{m} \Big)^{(2-\theta) / (1-\theta)} \ll K^{-\theta / (1-\theta)} Y_1.
\end{equation}

Next, we treat the non-diagonal part of $U_{1,2}$. We apply the first part of Lemma~\ref{Corput_kth_test} to the inner sum over $k$ (assuming it is non-empty) in~\eqref{U_1_2_split} with $\Lambda = RY_1 K^{-2}$ and $\nu = 2$. We get
\begin{multline*}
	\mathop{{\sum}^{'}}_{K \le k < 2K} e\bigl( c_3 k^{1/(1-\theta)} \eta \bigr) \ll \\
	K \bigl( RY_1 K^{-2} \bigr)^{1/2} + \bigl( RY_1 K^{-2} \bigr)^{-1/2} = (RY_1)^{1/2} + K (RY_1)^{-1/2}. 
\end{multline*} Thus, combining this bound with~\eqref{U_1_2_split} and~\eqref{diagonal}, we find
\begin{multline*}
	U_{1,2}(K) \ll K^{-\theta / (1-\theta)} Y_1 + \\
	\sum_R \Big( RKY_1^{\theta-1} (KY_1^{\theta-1})^{-(2-\theta)/(1-\theta)} \bigl[ (RY_1)^{1/2} + K (RY_1)^{-1/2} \bigr] \Big) \ll \\
	K^{-\theta / (1-\theta)} Y_1 + K^{1/2 - \theta / (1-\theta)} Y_1^{3\theta/2} + K^{-1/2 - \theta / (1-\theta)} Y_1^{\theta/2} \ll \\
	K^{-\theta / (1-\theta)} Y_1 + K^{1/2 - \theta / (1-\theta)} Y_1^{3\theta/2}.
\end{multline*} Hence, from~\eqref{S_1_2},
\begin{multline*}
	S_{1,2} (K) \ll \sum_{Y_1, Y_2} \Big( \sum_{K \le k < 2K} \max\Big( 1, \bigl(\log (kY_2^{\theta-1})\bigr)^2, \frac{Y_2^{2-\theta}}{k} \Big) k^{1/(1-\theta)} \Big)^{1/2} \cdot \\
	\bigl( K^{-\theta / (1-\theta)} Y_1 +  K^{1/2 - \theta / (1-\theta)} Y_1^{3\theta/2} \bigr)^{1/2} \ll \\
	\sum_{Y_1, Y_2} \bigl( K^{1/2} \log N + Y_2^{1-\theta/2} \bigr) \bigl( (KY_1)^{1/2} + (KY_1^{\theta})^{3/4} \bigr).
\end{multline*} Here we used the inequality $\max(1, (\log kY_2^{\theta-1})^2) \le (\log N)^2$. Finally,
\[
	S_{1,2} (K) \ll N^{1-\theta/2} \log N \bigl( (KN)^{1/2} + (KN^{\theta})^{3/4} \bigr),
\] and so
\begin{equation} \label{S_1_2_final_bound}
	\frac{2}{N^2} \sum_K S_{1,2} (K) \ll \frac{\log N}{N^2} \bigl( N^{2 - \theta / 2 + \varepsilon / 2} + N^{7/4 + \theta / 4 + 3\varepsilon / 4} \bigr),
\end{equation} which is $o(1)$ if $\varepsilon < \theta < 1 - 3\varepsilon$. The contribution from $S_{2,1} (K)$ is estimated similarly. \\

It remains only to evaluate the contribution from $S_{1,1}(K)$:
\[
	S_{1,1} (K) = \sum_{Y_1, Y_2} \sum_{K \le k < 2K} \hat{f} \Big( \frac{k}{N} \Big) T_k (Y_1) \overline{T_k (Y_2)} =: D_{1,1} (K) + E_{1,1} (K),
\] where $D_{1,1} (K)$ is the diagonal term ($Y_1 = Y_2$, $m_1 = m_2$),
\begin{multline*}
	D_{1,1} (K) = \\
	\sum_Y \sum_{K \le k < 2K} \hat{f} \Big( \frac{k}{N} \Big) \frac{1}{\alpha \theta (1-\theta)k} \bigl( \alpha k \theta \bigr)^{(2-\theta) / (1-\theta)} \sum_{\alpha \theta k (2Y)^{\theta - 1} \le m < \alpha \theta k Y^{\theta-1}} m^{-(2-\theta) / (1-\theta)}.
\end{multline*} By~\cite[Theorem~3.2]{Apostol}, the inner sum is
\[
	\frac{\bigl( \alpha \theta k Y^{\theta-1}\bigr)^{1 - (2-\theta) / (1-\theta)} - \bigl( 2^{\theta-1} \alpha \theta k Y^{\theta-1} \bigr)^{1 - (2-\theta)/(1-\theta)}}{1 - (2-\theta) / (1-\theta)} + O\bigl( (KY)^{-(2-\theta) / (1-\theta)} \bigr).
\] Simplifying the coefficients and using $\sum_Y Y = N + O(1)$, we obtain
\begin{multline*}
	S_{1,1} (K) = \sum_Y \sum_{K \le k < 2K} \hat{f} \Big( \frac{k}{N} \Big) Y + O\Big( \sum_Y \sum_{K \le k < 2K} k^{1 / (1-\theta)} (KY)^{-(2-\theta) / (1-\theta)} \Big) + \\ 
	E_{1,1} (K) = N \sum_{K \le k < 2K} \hat{f} \Big( \frac{k}{N} \Big) + O \bigl( K \bigr) + O \bigl( \log N \bigr) + E_{1,1} (K).
\end{multline*} Finally,
\begin{multline} \label{S_1_1_final_bound}
	\frac{2}{N^2} \sum_K S_{1,1} (K) = \\
	\frac{2}{N} \sum_{k \le N^{1+\varepsilon}} \hat{f} \Big( \frac{k}{N} \Big) + O\Big( \frac{1}{N} \Big) + O\Big( \frac{(\log N)^2}{N^2} \Big) + \frac{2}{N^2} \sum_K E_{1,1} (K) = \\
	\frac{1}{N} \sum_{k \in \Z} \hat{f} \Big( \frac{k}{N} \Big) + O_{\varepsilon} \Big( \frac{1}{N} \Big) + O\Big( \frac{(\log N)^2}{N^2} \Big) + \frac{2}{N^2} \sum_K E_{1,1} (K) = \\
	f(0) + O\Big( \frac{(\log N)^2}{N^2} \Big) + \frac{2}{N^2} \sum_K E_{1,1} (K).
\end{multline} To complete the proof of Proposition~\ref{prop_1}, we need to show that the sum with $E_{1,1} (K)$ is bounded by $S(N)$ given in~\eqref{new_expsum}.

\subsection{Poisson summation on $k$}

We have
\begin{multline*}
	E_{1,1} (K) = \sum_{Y_1, Y_2} \sum_{K \le k < 2K} \hat{f} \Big( \frac{k}{N} \Big) \sum_{\substack{m_1 \neq m_2 \\ \alpha \theta k (2Y_i)^{\theta-1} \le m_i < \alpha \theta k Y_i^{\theta-1}}} \Big(\frac{c_4}{\sqrt{k}} \Big)^2 \Big( \frac{k^2}{m_1 m_2} \Big)^{(2-\theta) / (2-2\theta)} \cdot \\
	e\Big( c_3 k^{1/(1-\theta)} \eta (m_1, m_2) \Big).
\end{multline*} We first remove the factor $\hat{f} (k/N)$ by partial summation:
\begin{lemma} \label{part_sum}
	Let $\{a_k\}$ and $\{ b_k \}$ be sequences of complex numbers, and $c > 1$ be a constant. If $T > 0$ is such that $|b_k - b_{k+1}| \le T/k$, then
	\[
		\Big| \sum_{A < k < B} a_k b_k \Big| \le \Big( \max_{A \le k \le cA} |b_k| + O(T) \Big) \max_{A \le \tilde A \le cA} \Big| \sum_{A \le k < \tilde A} a_k \Big|
	\] for any positive integers $A$ and $B$ satisfying $A \le B \le cA$. 
\end{lemma}

Now we apply Lemma~\ref{part_sum} to $E_{1,1} (K)$ with $b_k = \hat{f} (k/N)$ and
\[
	a_k = \sum_{\substack{m_1 \neq m_2 \\ \alpha \theta k (2Y_i)^{\theta-1} \le m_i < \alpha \theta k Y_i^{\theta-1}}} \Big(\frac{c_4}{\sqrt{k}} \Big)^2 \Big( \frac{k^2}{m_1 m_2} \Big)^{(2-\theta) / (2-2\theta)} e\Big( c_3 k^{1/(1-\theta)} \eta (m_1, m_2) \Big).
\] Since clearly $|\hat{f} (k/N) - \hat{f} ((k+1)/N)| \ll 1/N$, we can choose $T \asymp 1$. Thus,
\begin{multline*}
	E_{1,1} (K) \ll \sum_{Y_1, Y_2} \max_{K \le \tilde K \le 2K} \Big| \sum_{K \le k < \tilde K} c_4^2 k^{(2-\theta) / (1-\theta) - 1} \cdot \\
	\sum_{\substack{m_1 \neq m_2 \\ \alpha \theta k (2Y_i)^{\theta-1} \le m_i < \alpha \theta k Y_i^{\theta-1}}} \Big( \frac{1}{m_1 m_2} \Big)^{(2-\theta) / (2-2\theta)} e\Big( c_3 k^{1/(1-\theta)} \eta(m_1, m_2) \Big) \Big|.
\end{multline*} Next, we remove the factor $c_4^2 k^{(2-\theta) / (1-\theta) - 1}$ from the sum over $k$ by the second application of Lemma~\ref{part_sum} with

\begin{gather*}
	a_k = \sum_{\substack{m_1 \neq m_2 \\ \alpha \theta k (2Y_i)^{\theta-1} \le m_i < \alpha \theta k Y_i^{\theta-1}}} \Big( \frac{1}{m_1 m_2} \Big)^{(2-\theta) / (2-2\theta)} e\Big( c_3 k^{1/(1-\theta)} \eta(m_1, m_2) \Big), \\
	b_k = c_4^2 k^{(2-\theta) / (1-\theta) - 1}, \qquad T = c K^{(2-\theta) / (1-\theta) - 1}
\end{gather*} with some constant $c > 0$. Note that
\[
	\bigl| b_k - b_{k+1} \bigr| \ll \frac{1}{k} k^{(2-\theta) / (1-\theta) - 1}. 
\] Thus,
\begin{multline*}
	E_{1,1} (K) \ll K^{(2-\theta) / (1-\theta) - 1} \sum_{Y_1, Y_2} \max_{K \le \tilde K \le 2K} \Big| \sum_{K \le k < \tilde K} \\
	\sum_{\substack{m_1 \neq m_2 \\ \alpha \theta k (2Y_i)^{\theta-1} \le m_i < \alpha \theta k Y_i^{\theta-1}}} \Big( \frac{1}{m_1 m_2} \Big)^{(2-\theta) / (2-2\theta)}
	e\Big( c_3 k^{1/(1-\theta)} \eta(m_1, m_2) \Big) \Big|.
\end{multline*} Changing the order of summation, we obtain
\begin{multline} \label{SKKY1Y2}
	E_{1,1} (K) \ll K^{(2-\theta) / (1-\theta) - 1} \sum_{Y_1, Y_2} \max_{K \le \tilde K \le 2K} \\ 
	\Big| \sum_{\substack{m_1 \neq m_2 \\ \alpha \theta K (2Y_i)^{\theta-1} \le m_i < 2\alpha \theta K Y_i^{\theta-1}}} \Big( \frac{1}{m_1 m_2} \Big)^{(2-\theta)/(2-2\theta)} \sum_{K_1 \le k \le K_2} e\Big( c_3 k^{1/(1-\theta)} \eta(m_1, m_2) \Big) \Big|,
\end{multline} where $K_1$ and $K_2$ were defined in~\eqref{K_1} and~\eqref{K_2}.

Due to the symmetry between $m_1$ and $m_2$, we can add the restriction $m_1 < m_2$ to the corresponding sums in~\eqref{SKKY1Y2}. Then, applying Lemma~\ref{Poisson_summation} to the sum over $k$ in~\eqref{SKKY1Y2}, we get
\begin{multline} \label{E_1_1_bound}
	E_{1,1} (K) \ll K^{(2-\theta) / (1-\theta) - 1} \sum_{Y_1, Y_2} \max_{K \le \tilde K \le 2K} \Big| \sum_{\eta \in \mathcal{N}(K, Y_1, Y_2)} \Big( \frac{1}{m_1 m_2} \Big)^{(2-\theta) / (2-2\theta)} \cdot \\
	 \sum_{L \le \ell< \tilde L} \frac{c_1}{\sqrt{\eta}} \Big(\frac{\ell}{\eta}\Big)^{1/(2\theta) - 1} e\Big( c_2\ell^{1/\theta} \eta^{1 - 1/\theta} + \frac{1}{8} \Big) \Big|,
\end{multline} where $L, \tilde L, c_1, c_2$ were defined in Proposition~\ref{prop_1}. The factor $e(1/8)$ does not affect the absolute value of the inner double sum and, thus, can be omitted. We complete the proof by combining together the estimates~\eqref{S_2_2_final_bound}, \eqref{S_1_2_final_bound}, \eqref{S_1_1_final_bound} and~\eqref{E_1_1_bound}.

\section{Proof of Proposition~\ref{prop_2}}
\label{sec_pf_prop_2}

\subsection{Van der Corput differencing}

Let us denote the sum over $\eta \in \mathcal{N}(K, Y_1, Y_2)$ in~\eqref{E_1_1_bound} by $\tilde S(K, \tilde K, Y_1, Y_2)$, as in Proposition~\ref{prop_1}. Set the notation $m_2 =: m$, $m_1 =: m-r$. We then consider the sum over $r$. It has $\asymp KY_1^{\theta-1}$ terms with the restrictions $r \ll KY_2^{\theta-1}$ on the one side and (when $Y_1 > 2 Y_2$) $r \gg KY_2^{\theta-1}$ on the other side. We split it to the dyadic intervals $(R, 2R]$ as follows:
\[
	\bigl| \tilde S(K, \tilde K, Y_1, Y_2) \bigr| = \sum_{\substack{R \ll KY_2^{\theta-1} \\ \text{dyadic}}} \bigl| W(K, \tilde K, Y_1, Y_2, R) \bigr|,
\] where
\begin{multline*}
	W(K, \tilde K, Y_1, Y_2, R) := \\
	\sum_{\eta \in \mathcal{N}_R (K, Y_1, Y_2)} \Big( \frac{1}{(m-r)m} \Big)^{(2-\theta)/(2-2\theta)} \sum_{L \le \ell< \tilde L} \frac{c_1}{\sqrt{\eta}} \Big( \frac{\ell}{\eta} \Big)^{1/(2\theta) - 1} e\bigl( c_2\ell^{1/\theta} \eta^{1 - 1/\theta} \bigr)
\end{multline*} with
\begin{gather*}
	\mathcal{N}_R (K, Y_1, Y_2) := \bigl\{ \eta: \eta = m_1^{-\theta / (1-\theta)} - m_2^{-\theta / (1-\theta)}, \text{ where} \\
	\alpha \theta K (2Y_1)^{\theta-1} \le m_1 < 2 \alpha \theta K Y_1^{\theta-1}; \\  \alpha \theta K (2Y_2)^{\theta-1} \le m_2 < 2 \alpha \theta K Y_2^{\theta-1}; \\
	m_1 < m_2, \ m_2 - m_1 \in [R, 2R) \bigr\}.
\end{gather*}

Applying Cauchy inequality, we obtain
\begin{multline*}
	\bigl| W(K, \tilde K, Y_1, Y_2, R) \bigr|^2 \le \Big( \sum_{\eta \in \mathcal{N}_R (K, Y_1, Y_2)} \Big( \frac{1}{(m-r)m} \Big)^{(2-\theta) / (1-\theta)} \Big( \frac{1}{\eta} \Big)^{1 / \theta - 1} \Big) \Big) \cdot \\
	\Big( \sum_{\eta \in \mathcal{N}_R (K, Y_1, Y_2)} \Big| \sum_{L \le \ell< \tilde L}\ell^{1/(2\theta) - 1} e\bigl( c_2\ell^{1/\theta} \eta^{1 - 1/\theta} \bigr) \Big|^2 \Big).
\end{multline*} For the number of terms in $\mathcal{N}_R (K, Y_1, Y_2)$ we clearly have
\begin{equation} \label{bound_for_set_size}
	\# \mathcal{N}_R (K, Y_1, Y_2) \ll 
	\begin{cases}
		\displaystyle KY_1^{\theta-1} KY_2^{\theta-1},
		&\text{if } 
		\displaystyle Y_1 > 2Y_2;
		\\[2mm]
		\displaystyle R K Y_2^{\theta-1},
		&\text{if }
		\displaystyle Y_2 / 2 \le Y_1 \le 2Y_2.
	\end{cases}
\end{equation} In fact, $\# \mathcal{N}_R (K, Y_1, Y_2) = 0$ when $Y_1$ is much larger than $Y_2$, and, at the same time, $R$ is small compared to $KY_2^{\theta-1}$. Note also that $\# \mathcal{N}_R (K, Y_1, Y_2) = 0$ when $Y_1 \le Y_2 / 4$ by the definition. Next, recall that
\begin{equation} \label{size_of_L}
	L \asymp \eta K^{\theta / (1-\theta)}.
\end{equation} For $\eta \in \mathcal{N}_R (K, Y_1, Y_2)$ one has $\eta \asymp Z$, where 
\begin{equation} \label{size_of_Z}
	Z = Z(K, Y_1, Y_2, R) :=
	\begin{cases}
		\displaystyle Y_1^{\theta} K^{-\theta / (1-\theta)},
		&\text{if } 
		\displaystyle Y_1 > 2Y_2;
		\\[2mm]
		\displaystyle RY_2 K^{-1/(1-\theta)},
		&\text{if }
		\displaystyle Y_2 \le Y_1 \le 2Y_2.
	\end{cases}
\end{equation} Then, applying partial summation to the sum over $\ell$ (with $a_\ell= e(c_2\ell^{1/\theta} \eta^{1-1/\theta})$, $b_\ell=\ell^{1/(2\theta)-1}$, $T = L^{1/(2\theta)}$), we find
\begin{multline*}
	\bigl| W(K, \tilde K, Y_1, Y_2, R) \bigr|^2 \ll \\
	\Big( \#\mathcal{N}_R (K, Y_1, Y_2)\Big( \frac{1}{KY_1^{\theta-1} KY_2^{\theta-1}} \Big)^{(2-\theta) / (1-\theta)} \Big( \frac{1}{Z} \Big)^{1/\theta - 1} \Big) \cdot \\
	\bigl(ZK^{\theta / (1-\theta)}\bigr)^{1/\theta - 2} \Big( \sum_{\substack{\eta \in \mathcal{N}(K, Y_1, Y_2) \\ \eta \asymp Z}} \Big| \sum_{L \le \ell< L_1}  e\bigl( c_2\ell^{1/\theta} \eta^{1 - 1/\theta} \bigr) \Big|^2 \Big) = \\
	\# \mathcal{N}_R (K, Y_1, Y_2) K^{-3/(1-\theta)} (Y_1 Y_2)^{2-\theta} Z^{-1} \Big( \sum_{\substack{\eta \in \mathcal{N}(K, Y_1, Y_2) \\ \eta \asymp Z}} \Big| \sum_{L \le \ell< L_1}  e\bigl( c_2\ell^{1/\theta} \eta^{1 - 1/\theta} \bigr) \Big|^2 \Big).
\end{multline*} Then 
\begin{multline} \label{S_tilde_uniform}
	\bigl| \tilde S (K, \tilde K, Y_1, Y_2) \bigr| \ll \\
	\bigl(K^{- 3/(1-\theta)} (Y_1 Y_2)^{2-\theta} \bigr)^{1/2} \sum_{\substack{R \ll KY_2^{\theta-1} \\ \text{dyadic}}} \Big( \frac{\# \mathcal{N}_R (K, Y_1, Y_2)}{Z} \Big)^{1/2} U(K, \tilde K, Y_1, Y_2, R)^{1/2}, 
\end{multline} where
\begin{equation} \label{expression_for_U}
	U(K, \tilde K, Y_1, Y_2, R) := \sum_{\eta \in \mathcal{N}_R (K, Y_1, Y_2) } \Big| \sum_{L \le \ell< L_1}  e\bigl( c_2\ell^{1/\theta} \eta^{1 - 1/\theta} \bigr) \Big|^2.
\end{equation} 

To estimate $U(K, \tilde K, Y_1, Y_2, R)$, we first apply the van der Corput inequality (see~\cite[Lemma~8.17]{Iwaniec-Kowalski}):
\begin{lemma} \label{Corput_ineq}
	For any complex numbers $z_n$ we have
	\[
		\Big| \sum_{a < n < b} z_n \Big|^2 \le \Big( 1 + \frac{b-a}{H} \Big) \sum_{|h| < H} \Big( 1 - \frac{|h|}{H} \Big) \sum_{a < n, n+h < b} z_{n+h} \overline{z_n},
	\] where $H$ is any positive integer. 
\end{lemma}

Applying this lemma to the sum over $\ell$ in~\eqref{expression_for_U} with $H = H_1 = L^{\kappa_1}$ for some $0 < \kappa_1 < 1$, which will be chosen later, we obtain
\begin{multline*}
	\bigl| U(K, \tilde K, Y_1, Y_2, R) \bigr| \le \\
	\sum_{\eta \in \mathcal{N}_R (K, Y_1, Y_2) } \Big( 1 + \frac{L_1 - L}{H_1} \Big) \sum_{|h_1| < H_1} \Big( 1 - \frac{|h_1|}{H_1} \Big) \sum_{L \le \ell< L_1 - h_1} e\Big( c_2 \eta^{1-1/\theta} t(\ell, h_1) \Big),
\end{multline*} where 
\[
	t(\ell, h_1) := (\ell+h_1)^{1/\theta} -\ell^{1/\theta}.
\] Then
\begin{multline*}
	\bigl| U(K, \tilde K, Y_1, Y_2, R) \bigr| \ll \\ 
	\sum_{\eta \in \mathcal{N}_R (K, Y_1, Y_2) } \Big( \frac{L^2}{H_1} + \frac{L}{H_1} \sum_{0 < |h_1| < H_1} \Big| \sum_{L \le \ell< L_1 - h_1} e\Big( c_2 \eta^{1-1/\theta} t(\ell, h_1) \Big) \Big| \Big) \ll \\
	\# \mathcal{N}_R (K, Y_1, Y_2) \frac{L^2}{H_1} + \frac{L}{H_1} \sum_{\eta \in \mathcal{N}_R (K, Y_1, Y_2) } \sum_{0 < |h_1| < H_1} \Big| \sum_{L \le \ell< L_1 - h_1} e\Big( c_2 \eta^{1-1/\theta} t(\ell, h_1) \Big) \Big|.
\end{multline*} Applying Cauchy inequality again, we find
\begin{multline*}
	\bigl| U(K, \tilde K, Y_1, Y_2, R) \bigr| \le \# \mathcal{N}_R (K, Y_1, Y_2) \frac{L^2}{H_1} + 
	\frac{L}{H_1} \Big( \sum_{\eta \in \mathcal{N}_R (K, Y_1, Y_2) } \sum_{0 < |h_1| < H_1} 1 \Big)^{1/2} \cdot \\
	\Big( \sum_{\eta \in \mathcal{N}_R (K, Y_1, Y_2) } \sum_{0 < |h_1| < H_1} \Big| \sum_{L \le \ell< L_1 - h_1} e\Big( c_2 \eta^{1-1/\theta} t(\ell, h_1) \Big) \Big|^2 \Big)^{1/2}.
\end{multline*} Next, we apply Lemma~\ref{Corput_ineq} for the second time with $H = H_2 = L^{\kappa_2}$ for some $0 < \kappa_2 < 1$, which will be chosen later:
\begin{multline*}
	\bigl| U(K, \tilde K, Y_1, Y_2, R) \bigr| \ll \# \mathcal{N}_R (K, Y_1, Y_2) \frac{L^2}{H_1} + \frac{L}{H_1} \bigl( \#\mathcal{N}_R (K, Y_1, Y_2) H_1 \bigr)^{1/2} \cdot \\
	\Big( \sum_{\eta \in \mathcal{N}_R (K, Y_1, Y_2) } \sum_{0 \le |h_1| < H_1} \Big(1 + \frac{L_1 - h_1 - L}{H_2} \Big) \cdot \\
	\sum_{0 < |h_2| < H_2} \gamma(h_2) \sum_{L \le \ell< L_1 - h_1 - h_2} e\Big( c_2 \eta^{1 - 1/\theta} t(\ell, h_1, h_2) \Big) \Big)^{1/2},
\end{multline*} where $\gamma(h_2) := 1 - |h_2| / H_2$,
\[
	t(\ell, h_1, h_2) := \bigl( (\ell+h_2 + h_1)^{1/\theta} - (\ell+h_2)^{1/\theta} \bigr) - \bigl( (\ell+h_1)^{1/\theta} -\ell^{1/\theta} \bigr).
\] Hence,
\begin{multline*}
	\bigl| U(K, \tilde K, Y_1, Y_2, R) \bigr| \ll \\
	\# \mathcal{N}_R (K, Y_1, Y_2) \frac{L^2}{H_1} +
	\frac{L}{H_1} \bigl( \# \mathcal{N}_R (K, Y_1, Y_2) H_1 \bigr)^{1/2} \Big( \# \mathcal{N}_R (K, Y_1, Y_2) H_1\frac{L^2}{H_2} \Big)^{1/2} + \\ \frac{L}{H_1} \bigl( \# \mathcal{N}_R (K, Y_1, Y_2) H_1 \bigr)^{1/2} \Big( \frac{L}{H_2} \Big)^{1/2} \cdot \\
	\Big( \Big| \sum_{\substack{\eta, h_1, h_2 \\ \eta \in \mathcal{N}_R (K, Y_1, Y_2) \\ 0 < |h_i| < H_i}} \gamma(h_2) \sum_{L \le \ell< L_1 - h_1 - h_2} e\Big( c_2 \eta^{1 - 1/\theta} t(\ell, h_1, h_2) \Big) \Big| \Big)^{1/2}.
\end{multline*} We rewrite the last inequality as
\begin{multline} \label{bound_with_A}
	\bigl| U(K, \tilde K, Y_1, Y_2, R) \bigr| \le A_1(R) + A_2(R) + \\
	A_3(R) \Big( \Big| \sum_{\substack{\eta, h_1, h_2 \\ \eta \in \mathcal{N}_R (K, Y_1, Y_2) \\ 0 < |h_i| < H_i}} \gamma(h_2) \sum_{L \le \ell< L_1 - h_1 - h_2} e\Big( c_2 \eta^{1 - 1/\theta} t(\ell, h_1, h_2) \Big) \Big| \Big)^{1/2},
\end{multline} where
\begin{gather*}
	A_1(R) := \# \mathcal{N}_R (K, Y_1, Y_2) \frac{L^2}{H_1}, \qquad
	A_2(R) := \# \mathcal{N}_R (K, Y_1, Y_2) \frac{L^2}{H_2^{1/2}}, \\
	A_3(R) := \frac{L^{3/2} (\# \mathcal{N}_R (K, Y_1, Y_2))^{1/2}}{(H_1 H_2)^{1/2}}.
\end{gather*} Substituting the bound from~\eqref{bound_with_A}  into~\eqref{S_tilde_uniform}, we get
\begin{multline*}
	\bigl| \tilde S(K, \tilde K, Y_1, Y_2) \bigr| \ll \\
	\bigl(K^{- 3/(1-\theta)} (Y_1 Y_2)^{2-\theta} \bigr)^{1/2} \sum_R \Big( \frac{\# \mathcal{N}_R (K, Y_1, Y_2)}{Z} \Big)^{1/2} \Big( |A_1(R)|^{1/2} + \\
	|A_2(R)|^{1/2} + |A_3(R)|^{1/2} \Big( \Big| \sum_{\substack{\eta, h_1, h_2 \\ \eta \in \mathcal{N}_R (K, Y_1, Y_2)  \\ 0 < |h_i| < H_i}} \gamma(h_2) \sum_{L \le \ell< L_1 - h_1 - h_2} e\Big( c_2 \eta^{1 - 1/\theta} t(\ell, h_1, h_2) \Big) \Big| \Big)^{1/4} \Big).
\end{multline*} Then we have
\begin{multline*}
	\tilde S(K, \tilde K, Y_1, Y_2) \ll
	\sum_R \Big( B_1 (R) + B_2 (R) + \\
	B_3 (R) \Big( \Big| \sum_{\substack{\eta, h_1, h_2 \\ \eta \in \mathcal{N}_R (K, Y_1, Y_2) \\ 0 < |h_i| < H_i}} \gamma(h_2) \sum_{L \le \ell< L_1 - h_1 - h_2} e\Big( c_2 \eta^{1 - 1/\theta} t(\ell, h_1, h_2) \Big) \Big| \Big)^{1/4} \Big),
\end{multline*} where
\begin{gather}
	B_1 (R) := \Big( K^{-3 / (1-\theta)} (Y_1 Y_2)^{2-\theta}  \frac{\# \mathcal{N}_R^2 (K, Y_1, Y_2) L^2}{ZH_1}   \Big)^{1/2}, \label{B_1_R} \\
	B_2 (R) := \Big( K^{-3 / (1-\theta)} (Y_1 Y_2)^{2-\theta} \frac{\# \mathcal{N}_R^2 (K, Y_1, Y_2) L^2}{Z H_2^{1/2}} \Big)^{1/2}, \label{B_2_R} \\
	B_3 (R) := \Big( K^{-3 / (1-\theta)} (Y_1 Y_2)^{2-\theta} \frac{ \# \mathcal{N}_R^2 (K, Y_1, Y_2) L^{3/2}}{Z(H_1 H_2)^{1/2}}  \Big)^{1/2}. \label{B_3_R}
\end{gather} Upon inserting it into~\eqref{new_expsum}, we obtain
\begin{multline} \label{after_all_Corput}
	S(N) \ll \frac{1}{N^2} \sum_{K, Y_1, Y_2}  K^{(2-\theta)/(1-\theta) - 1} \max_{K \le \tilde K \le 2K} \sum_R \Big( B_1 (R) + B_2 (R) + \\
	B_3 (R) \Big( \Big| \sum_{\substack{\eta, h_1, h_2 \\ \eta \in \mathcal{N}_R (K, Y_1, Y_2) \\ 0 < |h_i| < H_i}} \gamma(h_2) \sum_{L \le \ell< L_1 - h_1 - h_2} e\Big( c_2 \eta^{1 - 1/\theta} t(\ell, h_1, h_2) \Big) \Big| \Big)^{1/4} \Big) =: \\
	S_{B_1} (N) + S_{B_2} (N) + S_{B_3} (N),
\end{multline} where the meaning of the notation $S_{B_i} (N)$ is clear. \\
	
Let us first estimate $S_{B_1} (N)$ and $S_{B_2} (N)$. \\

\textit{Case 1. $Y_2 /2  \le Y_1 \le 2 Y_2$.} Combining~\eqref{bound_for_set_size}, \eqref{size_of_L}, \eqref{size_of_Z}, and~\eqref{B_1_R}, we get
\begin{multline*}  
	\frac{1}{N^2} \sum_{\substack{K, Y_1, Y_2 \\ Y_2 / 2 \le Y_1 \le 2Y_2}} K^{(2-\theta)/(1-\theta)-1} \sum_R B_1 (R) \ll \\
	\frac{1}{N^2} \sum_{\substack{K, Y_1, Y_2 \\ Y_2 / 2 \le Y_1 \le 2Y_2}} K^{(2-\theta)/(1-\theta)-1} K^{-3/(2(1-\theta))} (Y_1 Y_2)^{1-\theta / 2} \cdot \\ \sum_R \frac{RKY_2^{\theta - 1} RY_2 K^{-1}}{\bigl( RY_2 K^{-1/(1-\theta)} \bigr)^{1/2} \bigl( RY_2 K^{-1} \bigr)^{\kappa_1/2} } \ll \frac{1}{N^2} \sum_{\substack{K, Y_1, Y_2 \\ Y_2 / 2 \le Y_1 \le 2Y_2}} K^{\kappa_1 / 2} Y_1^{1 - \theta / 2} Y_2^{1/2 + \theta / 2 - \kappa_1 / 2} \cdot \\
	\sum_R R^{3/2 - \kappa_1 / 2} \ll \frac{1}{N^2} \sum_{\substack{K, Y_1, Y_2 \\ Y_2 / 2 \le Y_1 \le 2Y_2}} K^{3/2} Y_1^{1 - \theta / 2} Y_2^{-1 + 2\theta - \kappa_1 \theta / 2} \ll N^{-1/2 + 3\theta / 2 - \kappa_1 \theta / 2 + 3\varepsilon / 2}.
\end{multline*}

Similarly,
\begin{multline*} 
	\frac{1}{N^2} \sum_{\substack{K, Y_1, Y_2 \\ Y_2 / 2 \le Y_1 \le 2Y_2}} K^{(2-\theta)/(1-\theta)-1} \sum_R B_2 (R) \ll \\
	\frac{1}{N^2} \sum_{\substack{K, Y_1, Y_2 \\ Y_2 / 2 \le Y_1 \le 2Y_2}} K^{(2-\theta)/(1-\theta)-1} K^{-3/(2(1-\theta))} (Y_1 Y_2)^{1-\theta / 2} \cdot \\ \sum_R \frac{RKY_2^{\theta-1} RY_2 K^{-1}}{\bigl( RY_2 K^{-1/(1-\theta)} \bigr)^{1/2} \bigl( RY_2 K^{-1} \bigr)^{\kappa_2 / 4}} \ll \frac{1}{N^2} \sum_{\substack{K, Y_1, Y_2 \\ Y_2 / 2 \le Y_1 \le 2Y_2}} K^{\kappa_2 / 4} Y_1^{1 - \theta / 2} Y_2^{1/2 - \kappa_2 / 4 + \theta / 2} \cdot \\
	\sum_R R^{3/2 - \kappa_2 / 4} \ll \frac{1}{N^2} \sum_{\substack{K, Y_1, Y_2 \\ Y_2 / 2 \le Y_1 \le 2Y_2}} K^{3/2} Y_1^{1 - \theta / 2} Y_2^{-1 + 2\theta - \kappa_2 \theta / 4} \ll N^{-1/2 + 3\theta /2 - \kappa_2 / 4 + 3\varepsilon / 2}.
\end{multline*}

\textit{Case 2. $Y_1 > 2Y_2$.} We have
\begin{multline*}
	\frac{1}{N^2} \sum_{\substack{K, Y_1, Y_2 \\ Y_1 > 2Y_2}} K^{(2-\theta)/(1-\theta)-1} \sum_R B_1 (R) \ll \\
	\frac{1}{N^2} \sum_{\substack{K, Y_1, Y_2 \\ Y_1 > 2Y_2}} K^{(2-\theta)/(1-\theta)-1} \sum_R K^{-3/(2(1-\theta))} (Y_1 Y_2)^{1-\theta/2} \frac{K^2 (Y_1 Y_2)^{\theta-1} Y_1^{\theta}}{Y_1^{\theta/2} K^{-\theta / (2(1-\theta))} Y_1^{\kappa_1 \theta/2}} \ll \\
	\frac{\log N}{N^2} \sum_{\substack{K, Y_1, Y_2 \\ Y_1 > 2Y_2}} K^{3/2} Y_1^{\theta - \kappa_1 \theta / 2} Y_2^{\theta / 2} \ll N^{-1/2 + 3\theta/2 - \kappa_1 \theta / 2 + 3\varepsilon / 2} \log N.
\end{multline*} Similarly,
\begin{multline*}
	\frac{1}{N^2} \sum_{\substack{K, Y_1, Y_2 \\ Y_1 > 2Y_2}} K^{(2-\theta)/(1-\theta)-1} \sum_R B_2 (R) \ll \\
	\frac{1}{N^2} \sum_{\substack{K, Y_1, Y_2 \\ Y_1 > 2Y_2}} K^{(2-\theta)/(1-\theta)-1} \sum_R K^{-3/(2(1-\theta))} (Y_1 Y_2)^{1-\theta/2} \frac{K^2 (Y_1 Y_2)^{\theta-1} Y_1^{\theta}}{Y_1^{\theta/2} K^{-\theta / (2(1-\theta))} Y_1^{\kappa_2 \theta/4}} \ll \\
	\frac{\log N}{N^2} \sum_{\substack{K, Y_1, Y_2 \\ Y_1 > 2Y_2}} K^{3/2} Y_1^{\theta / 2} Y_2^{\theta - \kappa_2 \theta / 4} \ll N^{-1/2 + 3\theta/2 - \kappa_2 \theta / 4 + 3\varepsilon / 2} \log N.
\end{multline*}

Thus,
\begin{gather} \label{S_B_1}
	S_{B_1} (N) \ll N^{-1/2 + 3\theta / 2 - \kappa_1 \theta / 2 + 3\varepsilon / 2} \log N, \\
 	\label{S_B_2}
	S_{B_2} (N) \ll N^{-1/2 + 3\theta / 2 - \kappa_2 \theta / 4 + 3\varepsilon / 2 } \log N.
\end{gather}

In the remainder of the paper, we will estimate $S_{B_3} (N)$.

\subsection{Double large sieve}

We handle the fivefold exponential sum over $\eta,\ell, h_1, h_2$ by combining the techniques of Bombieri-Iwaniec, Fourvy-Iwaniec-Cao-Zhai, and Robert-Sargos. Specifically, we use the double large sieve (as described in~\cite[Lemma~8]{Robert-Sargos}) to reduce the sum to two spacing inequalities.

\begin{lemma} \label{double_large_sieve}
	Define
	\begin{gather*}
		S = \sum_{k \le K} \sum_{\ell\le L} a(k) b(\ell) e\bigl( Xu(k)v(\ell) \bigr), \\
		\mathcal{B}_1 = \sum_{\substack{1 \le k_1, k_2 \le K \\ |u(k_1) - u(k_2)| \le X^{-1}}} 1, \qquad \mathcal{B}_2 = \sum_{\substack{1 \le\ell_1,\ell_2 \le L \\ |v(\ell_1) - v(\ell_2)| \le X^{-1}}} 1, 
	\end{gather*} where $|a(k)| \le 1$, $|b(\ell)| \le 1$ are complex valued functions, $|u(k)| \le 1$, $|v(\ell)| \le 1$ are real valued functions. Then
	\[
		S \ll X^{1/2} \mathcal{B}_1^{1/2} \mathcal{B}_2^{1/2}.
	\]
\end{lemma}

For convenience, we change the range of the summation over $h_1$ and $h_2$ in $S_{B_3}(N)$ in~\eqref{after_all_Corput} from $|h_i|<H_i$ to dyadic intervals $h_i\asymp \tilde{H}_i$, where $\tilde{H}_i$ denotes the size of the dyadic interval. There exist $\tilde H_1 \le H_1$ and $\tilde H_2 \le H_2$ such that
\begin{multline*}
	S_{B_3} (N) \ll \frac{\log N}{N^2} \sum_{K, Y_1, Y_2} K^{(2-\theta)/(1-\theta)-1} \max_{K \le \tilde K \le 2K} \sum_R B_3 (R) \cdot \\ \Big( \Big| \sum_{\substack{\eta, h_1, h_2 \\ \eta \in \mathcal{N}_R (K, Y_1, Y_2) \\ h_i \asymp \tilde H_i}} \gamma (h_2) \sum_{L \le \ell < L_1 - h_1 - h_2} e\Big( c_2 \eta^{1-1/\theta} t(\ell, h_1, h_2) \Big) \Big| \Big)^{1/4}.
\end{multline*}

Next, we apply Lemma~\ref{double_large_sieve} to the inner sum over $\eta, h_1, h_2,\ell$ in~\eqref{after_all_Corput} with $a \equiv 1$, $b = \gamma(h_2)$, $X = Z^{1 - 1/\theta} L^{1/\theta - 2} \tilde H_1 \tilde H_2$,
\begin{gather*}
	\mathcal{B}_1 = \sum_{\substack{\eta_1, \eta_2 \in \mathcal{N}_R (K, Y_1, Y_2) \\ |\eta_1^{\Theta} - \eta_2^{\Theta}| \le Z^{\Theta} X^{-1}}} 1, \qquad 
	\Theta := 1 - \frac{1}{\theta}, \\
	\mathcal{B}_2 = \sum_{\substack{L \le\ell, \tilde \ell< L_1 - h_1 - h_2 \\ h_1 \asymp \tilde H_1, \ h_2 \asymp \tilde  H_2 \\ |t(\ell, h_1, h_2) - t(\tilde\ell, \tilde h_1, \tilde h_2)| \le TX^{-1}}} 1, \qquad T = L^{1/\theta - 2} \tilde H_1 \tilde H_2.
\end{gather*} Then, using~\eqref{B_3_R}, we find
\begin{multline} \label{S_B_3}
	S_{B_3} (N)  \ll \frac{\log N}{N^2} \sum_{K, Y_1, Y_2}  K^{(2-\theta)/(1-\theta) - 1} \max_{K \le \tilde K \le 2K} \sum_R B_3 (R) X^{1/8} \mathcal{B}_1^{1/8} \mathcal{B}_2^{1/8} \ll \\
	\frac{1}{N^2} \sum_{K, Y_1, Y_2} K^{(2-\theta)/(1-\theta) - 1} \max_{K \le \tilde K \le 2K}  \sum_R K^{-3/(2(1-\theta))} (Y_1 Y_2)^{1 - \theta / 2} \cdot \\
	\frac{\bigl( \# \mathcal{N}_R (K, Y_1, Y_2) L \bigr)^{3/4}}{Z^{1/2} (H_1 H_2)^{1/4}} X^{1/8} \mathcal{B}_1^{1/8} \mathcal{B}_2^{1/8}.
\end{multline}

We estimate $\mathcal{B}_1$ and $\mathcal{B}_2$ separately. 
First, we replace $\mathcal{B}_1$ by a larger quantity. We drop the restriction $r \in [R, 2R)$ and set 
\[
	\eta_1 := m_1^{-\theta / (1-\theta)} - m_2^{-\theta / (1-\theta)}, \qquad \eta_2 := m_3^{-\theta / (1-\theta)} - m_4^{-\theta / (1-\theta)},
\] where $m_1, m_3 \asymp KY_1^{\theta-1}$ and $m_2, m_4 \asymp KY_2^{\theta-1}$. Applying the first-order Taylor approximation to 
\[
	\Big( \bigl( m_3^{-\theta / (1-\theta)} - m_4^{-\theta / (1-\theta)} \bigr)^{\Theta} + O(Z^{\Theta} X^{-1}) \Big)^{1/\Theta},
\] we get
\[
	m_1^{-\theta / (1-\theta)} - m_2^{-\theta / (1-\theta)} = m_3^{-\theta / (1-\theta)} - m_4^{-\theta / (1-\theta)} + O\bigl( ZX^{-1} \bigr).
\] Then
\[
	\mathcal{B}_1 \le \sum_{\substack{\alpha \theta K (2Y_1)^{\theta-1} \le m_1, m_3 < 2\alpha \theta KY_1^{\theta-1} \\ \alpha \theta K (2Y_2)^{\theta-1} \le m_2, m_4 < 2\alpha \theta KY_2^{\theta-1}  \\ |m_1^{\beta} - m_2^{\beta} - m_3^{\beta} + m_4^{\beta}| \le c Z X^{-1}}} 1 \le \sum_{\substack{1 \le m_1, \ldots, m_4 < 4\alpha \theta KY_2^{\theta-1} \\ |m_1^{\beta} - m_2^{\beta} - m_3^{\beta} + m_4^{\beta}| \le c Z X^{-1}}} 1,
\] where $\beta := -\theta / (1-\theta)$ and $c > 0$ is the appropriate constant. To estimate this quantity, we use a slight variation of Robert and Sargos's inequality (see~\cite[Theorem~2]{Robert-Sargos}):
\begin{lemma} \label{Robert-Sargos}
	Let $M \ge 2$ be an integer, $a_2 > a_1 \ge 0$, $\delta > 0$ and $\alpha \ne 0, 1$ be real numbers, $\mathcal{N}(M, \delta)$ be the number of quadruplets $(m_1, m_2, m_3, m_4)$ of integers, $a_1 M \le m_i < a_2M$, $i = 1, \ldots, 4$, satisfying the inequality
	\[
		|m_1^{\alpha} + m_2^{\alpha} - m_3^{\alpha} - m_4^{\alpha}| \le \delta M^{\alpha}.
	\] Then
	\[
		\mathcal{N}(M, \delta) \ll_{\varepsilon} M^{2+\varepsilon} + \delta M^{4+\varepsilon}.
	\]
\end{lemma} Thus,
\begin{equation} \label{B_1_spacing_bound}
	\mathcal{B}_1 \ll_{\varepsilon_1} (KY_2^{\theta-1})^{2+\varepsilon_1} + \bigl((KY_2^{\theta-1})^{\theta / (1-\theta)} ZX^{-1}\bigr) (KY_2^{\theta-1})^{4+\varepsilon_1}
\end{equation} with some $\varepsilon_1 > 0$.

We estimate $\mathcal{B}_2$ using the result of Cao and Zhai (see~\cite[Theorems~1 and~2]{Cao-Zhai}):
\begin{lemma} \label{Cao-Zhai}
	Let $M, H_1, H_2$ be integers $\ge 10$ such that $H_1 \le H_2 \le M^{2/3 - \varepsilon}$, $\Delta > 0$, $\alpha \neq 0,1,2,3$ be real, and let $T := M^{\alpha-2} H_1 H_2$. Furthermore, let $\mathcal{F}(M, H_1, H_2, \Delta)$ be the number of sextuplets $(m, \tilde m, h_1, \tilde h_1, h_2, \tilde h_2)$ with $M \le m, \tilde m < 2M$, $h_1, \tilde h_1 \asymp H_1$, and $h_2, \tilde h_2 \asymp H_2$, satisfying
	\[
		\bigl|t(m, h_1, h_2) - t(\tilde m, \tilde h_1, \tilde h_2) \bigr| \le \Delta T.
	\] Then
	\[
		\mathcal{F}(M, H_1, H_2, \Delta) M^{-2 \varepsilon} \ll MH_1 H_2 + \Delta (M H_1 H_2)^2 + M^{-2} H_1^2 H_2^6 + H_1^2 H_2^{8/3}.
	\] Alternatively, if $H_1$ is much smaller than $H_2$, namely one has $H_1 M^{\varepsilon} < H_2 \le M^{1-\varepsilon}$, and $H_1 H_2 \le M^{3/2 - \varepsilon}$, then one also has the inequality
	\begin{multline*}
		\mathcal{F}(M, H_1, H_2, \Delta) M^{-2 \varepsilon} \ll MH_1 H_2 + \Delta (M H_1 H_2)^2 + (M H_1^7 H_2^9)^{1/4} + M^{-2} H_1^4 H_2^4 + \\
		(\Delta M^4 H_1^{15} H_2^{17})^{1/8} + (\Delta M^4 H_1^3 H_2)^{1/2} + (H_1^{13} H_2^{15})^{1/6} + (\Delta M^2 H_1^8 H_2^{10})^{1/4} + (M^{-1} H_1^5 H_2^6)^{1/2}.
	\end{multline*}
\end{lemma} We first assume that $\tilde H_1$ and $\tilde H_2$ are separated at least by the factor of $L^{\varepsilon_1}$. Precisely, we have either $\tilde H_1 L^{\varepsilon_1} \ll \tilde H_2$ or $\tilde H_2 L^{\varepsilon_1} \ll \tilde H_1$, where $\varepsilon_1$ will be chosen later. The case $\tilde H_1 \approx \tilde H_2$ will be considered at the end of the paper. Applying the second part of Lemma~\ref{Cao-Zhai} to $\mathcal{B}_2$ we obtain
\begin{multline*}
	\mathcal{B}_2 \ll_{\varepsilon_1} L^{2\varepsilon_1} \Big( L \tilde H_1 \tilde H_2 +  X^{-1} (L \tilde H_1 \tilde H_2)^2 + (L \tilde H_1^7 \tilde H_2^9)^{1/4} + L^{-2} \tilde H_1^4 \tilde H_2^4 + (X^{-1} L^4 \tilde H_1^{15} \tilde H_2^{17})^{1/8} + \\ (X^{-1} L^4 \tilde H_1^3 \tilde H_2)^{1/2} + (\tilde H_1^{13} \tilde H_2^{15})^{1/6} + (X^{-1} L^2 \tilde H_1^8 \tilde H_2^{10})^{1/4} + (L^{-1} \tilde H_1^5 \tilde H_2^6)^{1/2} \Big).
\end{multline*} Here we make the same choice of $\varepsilon_1$. 

Next, note that the expression $X^{1/8} \mathcal{B}_1^{1/8} \mathcal{B}_2^{1/8}$ is maximized when $\tilde H_1$ and $\tilde H_2$ take their maximum values, that is, when $\tilde H_1 = H_1$ and $\tilde H_2 = H_2$. Therefore,
\begin{multline*}
	S_{B_3} (N) \ll \frac{N^{\varepsilon + \varepsilon_1}}{N^2} \sum_{K, Y_1, Y_2} K^{(2-\theta) / (1-\theta) - 1} \max_{K \le \tilde K \le 2K} \sum_R K^{-3/(2(1-\theta))} (Y_1 Y_2)^{1-\theta/2} \cdot \\
	\frac{\bigl( \# \mathcal{N}_R (K, Y_1, Y_2) L \bigr)^{1/4}}{Z^{1/2} (H_1 H_2)^{1/4}} \bigl( Z^{1-1/\theta} L^{1/\theta - 2} H_1 H_2 \bigr)^{1/8} \cdot \\
	\Big( (KY_2^{\theta-1})^2 + \bigl( (KY_2^{\theta-1})^{\theta / (1-\theta)} ZX^{-1} \bigr) (KY_2^{\theta-1})^4 \Big)^{1/8} \cdot \\
	 \Big( L H_1 H_2 +  X^{-1} (L H_1 H_2)^2 + (L H_1^7 H_2^9)^{1/4} + L^{-2} H_1^4 H_2^4 + (X^{-1} L^4 H_1^{15} H_2^{17})^{1/8} + \\ (X^{-1} L^4 H_1^3 H_2)^{1/2} + ( H_1^{13} H_2^{15})^{1/6} + (X^{-1} L^2 H_1^8 H_2^{10})^{1/4} + (L^{-1} H_1^5 H_2^6)^{1/2} \Big)^{1/8}.
\end{multline*}

Next, we consider two cases once again. \\

\textit{Case 1. $Y_2 / 2 \le Y_1 \le 2Y_2$.} \\

From ~\eqref{bound_for_set_size}, \eqref{size_of_L}, and~\eqref{size_of_Z} we have 
\begin{gather*}
	\# \mathcal{N}_R (K, Y_1, Y_2) \ll RK Y_2^{\theta-1}, \qquad Z \asymp R Y_2 K^{-1/(1-\theta)}, \qquad   L \asymp ZK^{\theta / (1-\theta)} \asymp R Y_2 K^{-1}, \\
	X = Z^{1-1/\theta} L^{1/\theta - 2} H_1 H_2 = K^2 Y_2^{-1} R^{-1} H_1 H_2.
\end{gather*} Then 
\begin{multline} \label{S_B_3_1_before_R}
	S_{B_3}^{(1)} (N) := S_{B_3} (N; \ Y_2 / 2 \le Y_1 \le 2 Y_2) \ll \\
	\frac{N^{\varepsilon + \varepsilon_1}}{N^2} \sum_{\substack{K, Y_1, Y_2 \\ Y_2 / 2 \le Y_1 \le 2Y_2}} Y_1^{1 - \theta / 2} Y_2^{1/2 + \theta / 4} \sum_R \frac{R}{(H_1 H_2)^{1/4}} \cdot \\
	\bigl( K^2 Y_2^{-1} R^{-1} H_1 H_2 \bigr)^{1/8} \Big( (KY_2^{\theta-1})^2 + \frac{(KY_2^{\theta-1})^4}{K^3 Y_2^{\theta - 2} R^{-2} H_1 H_2} \Big)^{1/8} \cdot \\
	\Big( RY_2K^{-1} H_1 H_2 + \frac{(RY_2 K^{-1} H_1 H_2)^2}{K^2 Y_2^{-1} R^{-1} H_1 H_2} + (RY_2 K^{-1} H_1^7 H_2^9)^{1/4} + (RY_2 K^{-1})^{-2} H_1^4 H_2^4 + \\
	\Big(\frac{ (RY_2 K^{-1})^4 H_1^{15} H_2^{17}}{K^2 Y_2^{-1} R^{-1} H_1 H_2} \Big)^{1/8} + \Big( \frac{(RY_2 K^{-1})^4 H_1^3 H_2}{K^2 Y_2^{-1} R^{-1} H_1 H_2} \Big)^{1/2} + ( H_1^{13} H_2^{15})^{1/6} + \\
	\Big(\frac{(RY_2 K^{-1})^2 H_1^8 H_2^{10}}{K^2 Y_2^{-1} R^{-1} H_1 H_2}\Big)^{1/4} + ((RY_2 K^{-1})^{-1} H_1^5 H_2^6)^{1/2} \Big)^{1/8}.
\end{multline}

The most contribution to the sum over $R$ comes from $R$ maximal, namely $R \asymp KY_2^{\theta-1}$. It corresponds to $Z \asymp Y_2^{\theta} K^{-\theta / (1-\theta)}$, $X \asymp KY_2^{-\theta} H_1 H_2$, $L \asymp Y_2^{\theta}$, $H_1 \asymp Y_2^{\kappa_1 \theta}$, $H_2 \asymp Y_2^{\kappa_2 \theta}$. Then, we get
\begin{multline*}
	S_{B_3}^{(1)} (N) \ll \frac{N^{\varepsilon + \varepsilon_1}}{N^2} \sum_{\substack{K, Y_1, Y_2 \\ Y_2 /2 \le Y_1 \le 2Y_2}} \frac{KY_2^{1/2 + 3\theta/4}}{(H_1 H_2)^{1/4}} \cdot \\
	\bigl( KY_2^{-\theta} H_1 H_2 \bigr)^{1/8} \Big( (KY_2^{\theta-1})^2 + \frac{(KY_2^{\theta-1})^4}{KY_2^{-\theta} H_1 H_2} \Big)^{1/8} \cdot \\
	\Big( Y_2^{\theta} H_1 H_2 + \frac{(Y_2^{\theta} H_1 H_2)^2}{KY_2^{-\theta} H_1 H_2} + (Y_2^{\theta} H_1^7 H_2^9)^{1/4} + Y_2^{-2\theta} H_1^4 H_2^4 + \\
	\Big(\frac{ Y_2^{4\theta} H_1^{15} H_2^{17}}{KY_2^{-\theta} H_1 H_2} \Big)^{1/8} + \Big( \frac{Y_2^{4\theta} H_1^3 H_2}{KY_2^{-\theta} H_1 H_2} \Big)^{1/2} + ( H_1^{13} H_2^{15})^{1/6} + \\
	\Big(\frac{Y_2^{2\theta} H_1^8 H_2^{10}}{KY_2^{-\theta} H_1 H_2}\Big)^{1/4} + ( Y_2^{-\theta} H_1^5 H_2^6)^{1/2} \Big)^{1/8}.
\end{multline*}	

The maximum contribution to the sum over $K, Y_1, Y_2$ comes from $Y_2$ (which is $\asymp Y_1$) and $K$ maximal, namely $Y_1, Y_2 \asymp N$, $K \asymp N^{1+\varepsilon}$. Thus,
\begin{multline} \label{S_B_3_1}
	S_{B_3}^{(1)} (N) \ll N^{3\varepsilon + \varepsilon_1} \frac{N^{-1/2 + 3\theta / 4}}{(H_1 H_2)^{1/4}}  \bigl( N^{1-\theta} H_1 H_2 \bigr)^{1/8} \Big( N^{2\theta} + \frac{N^{4\theta}}{N^{1-\theta} H_1 H_2} \Big)^{1/8} \cdot \\
	\Big( N^{\theta} H_1 H_2 + \frac{\bigl( N^{\theta} H_1 H_2 \bigr)^2}{N^{1-\theta} H_1 H_2 } + (N^{\theta} H_1^7 H_2^9)^{1/4} + N^{-2\theta} H_1^4 H_2^4 + ( N^{5\theta-1} H_1^{14} H_2^{16} )^{1/8} + \\ 
	( N^{5\theta - 1} H_1^2 )^{1/2} + ( H_1^{13} H_2^{15})^{1/6} + ( N^{3\theta - 1} H_1^7 H_2^9 )^{1/4} + ( N^{-\theta} H_1^5 H_2^6)^{1/2} \Big)^{1/8}.
\end{multline}

\textit{Case 2. $Y_1 > 2Y_2$.} \\

From ~\eqref{bound_for_set_size}, \eqref{size_of_L}, and~\eqref{size_of_Z} we have
\begin{gather*}
	\# \mathcal{N}_R (K, Y_1, Y_2) \ll KY_1^{\theta-1} KY_2^{\theta-1}, \qquad Z = Y_1^{\theta} K^{-\theta / (1-\theta)}, \qquad L \asymp Z K^{\theta / (1-\theta)} = Y_1^{\theta},  \\
	X \asymp KY_1^{-\theta} H_1 H_2.
\end{gather*} Simplifying the right hand side in~\eqref{S_B_3}, we obtain
\begin{multline*}
	S_{B_3}^{(2)} (N) := S_{B_3} \bigl( N; \ Y_1 > 2Y_2 \bigr) \ll \\
	\frac{\log N}{N^2} \sum_{\substack{K, Y_1, Y_2 \\ Y_1 > 2Y_2}} \frac{K Y_1^{1/4 + \theta/2} Y_2^{1/4 + \theta / 4}}{(H_1 H_2)^{1/4}} \bigl( KY_1^{-\theta} H_1 H_2 \bigr)^{1/8} \mathcal{B}_1^{1/8} \mathcal{B}_2^{1/8},
\end{multline*} where ($\log N$)-factor comes from the sum over $R$. 

Here we need to consider two more cases. First, we assume that the first term in the bound~\eqref{B_1_spacing_bound} for $\mathcal{B}_1$ is dominating: \\

\textit{Case 2.1. $(KY_2^{\theta-1})^{2+\varepsilon_1} \gg \bigl( (KY_2^{\theta-1})^{\theta / (1-\theta)} ZX^{-1} \bigr) (KY_2^{\theta-1})^{4+\varepsilon_1}$.} \\

This condition can be rewritten as	$Y_2^{2-\theta} \gg KY_1^{2\theta} (H_1 H_2)^{-1}$. We have
\begin{multline} \label{S_B_3_2_1_before_KY}
	S_{B_3}^{(2.1)} (N) := S_{B_3} \bigl( N; \ Y_1 > 2Y_2; \ Y_2^{2-\theta} \gg KY_1^{2\theta} (H_1 H_2)^{-1} \bigr) \ll \\
	\frac{N^{\varepsilon + \varepsilon_1}}{N^2} \sum_{\substack{K, Y_1, Y_2 \\ Y_1 > 2Y_2 \\ Y_2^{2-\theta} \gg KY_1^{2\theta} (H_1 H_2)^{-1}}} \frac{K Y_1^{1/4 + \theta/2} Y_2^{1/4 + \theta / 4}}{(H_1 H_2)^{1/4}} \bigl( KY_1^{-\theta} H_1 H_2 \bigr)^{1/8} \cdot \\ 
	\bigl( KY_2^{\theta-1} \bigr)^{1/4} \Big( Y_1^{\theta} H_1 H_2 + \frac{(Y_1^{\theta} H_1 H_2)^2}{KY_1^{-\theta} H_1 H_2} + (Y_1^{\theta} H_1^7 H_2^9)^{1/4} + Y_1^{-2\theta} H_1^4 H_2^4 + \\
	\Big(\frac{ Y_1^{4\theta} H_1^{15} H_2^{17}}{KY_1^{-\theta} H_1 H_2} \Big)^{1/8} + \Big( \frac{Y_1^{4\theta} H_1^3 H_2}{KY_1^{-\theta} H_1 H_2} \Big)^{1/2} + ( H_1^{13} H_2^{15})^{1/6} + \\
	\Big(\frac{Y_1^{2\theta} H_1^8 H_2^{10}}{KY_1^{-\theta} H_1 H_2}\Big)^{1/4} + ( Y_1^{-\theta} H_1^5 H_2^6)^{1/2} \Big)^{1/8}.
\end{multline} Recall that $H_1 H_2 = Y_1^{(\kappa_1 + \kappa_2) \theta}$ with $\kappa_1 + \kappa_2 \le \frac{3}{2} - \varepsilon_1$. The most contribution in the last upper bound for $S_{B_3}^{(2.1)} (N)$ comes from $K, Y_1, Y_2$ maximal. We get
\begin{multline} \label{S_B_3_2_1}
	S_{B_3}^{(2.1)} (N) \ll N^{3\varepsilon + \varepsilon_1} \frac{N^{-1/2 + 3\theta / 4}}{(H_1 H_2)^{1/4}} \bigl( N^{1-\theta} H_1 H_2 \bigr)^{1/8} N^{\theta / 4} \cdot \\
		\Big( N^{\theta} H_1 H_2 + \frac{\bigl( N^{\theta} H_1 H_2 \bigr)^2}{N^{1-\theta} H_1 H_2 } + (N^{\theta} H_1^7 H_2^9)^{1/4} + N^{-2\theta} H_1^4 H_2^4 + ( N^{5\theta-1} H_1^{14} H_2^{16} )^{1/8} + \\ 
	( N^{5\theta - 1} H_1^2 )^{1/2} + ( H_1^{13} H_2^{15})^{1/6} + ( N^{3\theta - 1} H_1^7 H_2^9 )^{1/4} + ( N^{-\theta} H_1^5 H_2^6)^{1/2} \Big)^{1/8}.
\end{multline}

\textit{Case 2.2. $(KY_2^{\theta-1})^{2+\varepsilon_1} \ll \bigl( (KY_2^{\theta-1})^{\theta / (1-\theta)} ZX^{-1} \bigr) (KY_2^{\theta-1})^{4+\varepsilon_1}$.} \\

Similarly, this condition can be rewritten as $Y_2^{2-\theta} \ll  KY_1^{2\theta} (H_1 H_2)^{-1}$. Since the second term in the bound for $\mathcal{B}_1$~\eqref{B_1_spacing_bound} dominates in this case, we obtain:
\begin{multline} \label{S_B_3_2_2_first}
	S_{B_3}^{(2.2)} (N) := S_{B_3} \bigl( N; \ Y_1 > 2Y_2; \ Y_2^{2-\theta} \ll KY_1^{2\theta} (H_1 H_2)^{-1} \bigr) \ll \\
	\frac{N^{\varepsilon + \varepsilon_1}}{N^2} \sum_{\substack{K, Y_1, Y_2 \\ Y_1 > 2Y_2 \\ Y_2^{2-\theta} \ll KY_1^{2\theta} (H_1 H_2)^{-1}}} \frac{K Y_1^{1/4 + \theta/2} Y_2^{1/4 + \theta / 4}}{(H_1 H_2)^{1/4}} \bigl( KY_1^{-\theta} H_1 H_2 \bigr)^{1/8} \cdot \\ 
	\Big( \frac{(KY_2^{\theta-1})^4}{KY_1^{-2\theta}Y_2^{\theta} H_1 H_2} \Big)^{1/8} \Big( Y_1^{\theta} H_1 H_2 + \frac{(Y_1^{\theta} H_1 H_2)^2}{KY_1^{-\theta} H_1 H_2} + (Y_1^{\theta} H_1^7 H_2^9)^{1/4} + Y_1^{-2\theta} H_1^4 H_2^4 + \\
	\Big(\frac{ Y_1^{4\theta} H_1^{15} H_2^{17}}{KY_1^{-\theta} H_1 H_2} \Big)^{1/8} + \Big( \frac{Y_1^{4\theta} H_1^3 H_2}{KY_1^{-\theta} H_1 H_2} \Big)^{1/2} + ( H_1^{13} H_2^{15})^{1/6} + \\
	\Big(\frac{Y_1^{2\theta} H_1^8 H_2^{10}}{KY_1^{-\theta} H_1 H_2}\Big)^{1/4} + ( Y_1^{-\theta} H_1^5 H_2^6)^{1/2} \Big)^{1/8}.
\end{multline}

The maximum contribution to this expression comes from $K, Y_1$ maximal and $Y_2 \asymp 1$. We get
\begin{multline} \label{S_B_3_2_2_final}
	S_{B_3}^{(2.2)} (N) \ll N^{3\varepsilon + \varepsilon_1} \frac{N^{-1/4 + 5\theta / 8}}{(H_1 H_2)^{1/4}} \Big(N^{\theta} H_1 H_2 + N^{3\theta-1} H_1 H_2 + (N^{\theta} H_1^7 H_2^9)^{1/4} + \\ 
	 N^{-2\theta} H_1^4 H_2^4 + ( N^{5\theta-1} H_1^{14} H_2^{16} )^{1/8} + 
	( N^{5\theta - 1} H_1^2 )^{1/2} + ( H_1^{13} H_2^{15})^{1/6} + \\
	( N^{3\theta - 1} H_1^7 H_2^9 )^{1/4} + ( N^{-\theta} H_1^5 H_2^6)^{1/2} \Big)^{1/8}.
\end{multline} \\

Combining~\eqref{S_B_3_1} and~\eqref{S_B_3_2_1}, we find
\begin{multline} \label{S_B_3_1_and_2_1}
	S_{B_3}^{(1)} (N) + S_{B_3}^{(2.1)} (N) \ll \\
	N^{3\varepsilon + \varepsilon_1} \frac{N^{-1/2 + 3\theta / 4}}{(H_1 H_2)^{1/4}} \bigl( N^{1-\theta} H_1 H_2 \bigr)^{1/8} \Big( N^{2\theta} + \frac{N^{4\theta}}{N^{1-\theta} H_1 H_2} \Big)^{1/8} \cdot \\
	\Big( N^{\theta} H_1 H_2 + N^{3\theta-1} H_1 H_2 + (N^{\theta} H_1^7 H_2^9)^{1/4} + N^{-2\theta} H_1^4 H_2^4 + ( N^{5\theta-1} H_1^{14} H_2^{16} )^{1/8} + \\ 
	( N^{5\theta - 1} H_1^2 )^{1/2} + ( H_1^{13} H_2^{15})^{1/6} + ( N^{3\theta - 1} H_1^7 H_2^9 )^{1/4} + ( N^{-\theta} H_1^5 H_2^6)^{1/2} \Big)^{1/8}.
\end{multline}

We can drop the term $N^{4\theta} / (N^{1-\theta} H_1 H_2)$ in the second factor by assuming that $H_1 H_2 \gg N^{3\theta-1}$. Choose $\varepsilon_1 < \varepsilon$. Then, from~\eqref{after_all_Corput} and~\eqref{S_B_3_1_and_2_1},
\[
	S(N) \ll S_{B_1} (N) + S_{B_2} (N) + N^{4\varepsilon} \sum_{i=1}^9 S_i + S_{B_3}^{(2.2)} (N),
\] where $S_i$ correspond to nine terms in the bound in~\eqref{S_B_3_1_and_2_1}. We have
\[
	S_1 = N^{-1/2 + 3\theta/4 + 1/8 + \theta/4}, \qquad S_2 = N^{-1/2 + 5\theta/4}.
\] Then $S_1, S_2 = o(N^{-4\varepsilon})$ if $\theta < \frac{3}{8} - 4\varepsilon = 0.375 - 4\varepsilon$. Similarly,
\begin{gather*}
	S_3 = N^{-3/8 + 29\theta / 32} H_1^{3/32} H_2^{5/32}, \qquad
	S_4 = N^{-3/8 + 5\theta / 8} H_1^{3/8} H_2^{3/8}, \\
	S_5 = N^{-25/64 + 61\theta / 64} H_1^{6/64} H_2^{8/64}, \qquad
	S_6 = N^{-7/16 + 19\theta / 16} H_2^{-1/8}, \\
	S_7 = N^{-3/8 + 7\theta / 8} H_1^{13/48} H_2^{15/48}, \qquad
	S_8 = N^{-13/32 + 31\theta / 32} H_1^{3/32} H_2^{5/32}, \\
	S_9 = N^{-3/8 + 13\theta / 16} H_1^{3/16} H_2^{4/16}.
\end{gather*}

By combining the expressions with the bounds~\eqref{S_B_1} and~\eqref{S_B_2} for $S_{B_1} (N)$ and $S_{B_2} (N)$, we obtain that $S(N) = o(1)$ when
\begin{gather*}
	\theta < \min \Big( \frac{1}{3-\kappa_1}, \ \frac{2}{6 - \kappa_2}, \ \frac{12}{29+3\kappa_1+5\kappa_2}, \ 
	\frac{3}{5 + 3\kappa_1 + 3\kappa_2}, \
	\frac{25}{61+6\kappa_1+8\kappa_2}, \
	\frac{7}{19 - 2\kappa_2}, \\
	\frac{18}{42+7\kappa_1+9\kappa_2}, \ 
	\frac{13}{31+3\kappa_1+5\kappa_2}, \
	\frac{6}{13+3\kappa_1+4\kappa_2} \Big) - 5 \varepsilon,
\end{gather*} which is maximized at $(\kappa_1, \kappa_2) = (\frac{12}{43}, \frac{24}{43})$. Therefore, we have $\theta < \frac{43}{117} - 5\varepsilon$, with the largest term being $S_7$. From~\eqref{S_B_3_2_2_final} we find
\begin{multline*}
	S_{B_3}^{(2.2)} (N) \ll N^{4\varepsilon} \Big( 
	N^{-1/4 + 3\theta/4 - (\kappa_1 + \kappa_2)\theta / 8} +
	N^{-3/8 + \theta - (\kappa_1 + \kappa_2)\theta / 8} + \\
	N^{-1/4 + 21\theta/32 - \kappa_1 \theta / 32 + \kappa_2 \theta / 32} +
	N^{-1/4 + 3\theta/8 + (\kappa_1 + \kappa_2)\theta / 4} +
	N^{-17/64 + 45\theta/64 - 16\kappa_1 \theta / 64} + \\
	N^{-5/16 + 15\theta/16 - 2\kappa_1 \theta / 16 - 4\kappa_2 \theta / 16} +
	N^{-1/4 + 5\theta/8 + \kappa_1 \theta / 48 + 3\kappa_2 \theta / 48} + \\
	N^{-9/32 + 23\theta/32 - \kappa_1 \theta / 32 + \kappa_2 \theta / 32} + 
	N^{-1/4 + 9\theta/16 + \kappa_1 \theta / 16 + 2\kappa_2 \theta / 16} \Big),
\end{multline*} which is $o(1)$ when
\begin{gather*}
	\theta < \min \Big( 
	\frac{2}{6-\kappa_1 - \kappa_2}, \ 
	\frac{3}{8 - \kappa_1 - \kappa_2}, \ 
	\frac{8}{21 - \kappa_1 + \kappa_2}, \ 
	\frac{2}{3 + 2\kappa_1 + 2\kappa_2}, \
	\frac{17}{45 - 16\kappa_1}, \\
	\frac{5}{15 - 2\kappa_1 - 4\kappa_2}, \ 
	\frac{12}{30 + \kappa_1 + 3\kappa_2}, \ 
	\frac{9}{23 - \kappa_1 + \kappa_2}, \
	\frac{4}{9 + \kappa_1 + 2\kappa_2} \Big) - 5 \varepsilon.
\end{gather*} It covers the interval $\theta \in [0, \frac{43}{117})$ when $(\kappa_1, \kappa_2) = (\frac{12}{43}, \frac{24}{43})$. Note that with this choice of $(\kappa_1, \kappa_2)$ the condition $H_1 H_2 \gg N^{3\theta-1}$ is also satisfied.

\subsection{Case $\tilde H_1 \approx \tilde H_2$} We assume that $L^{-\varepsilon_1} \tilde H_1 \ll \tilde H_2 \ll L^{\varepsilon_1} \tilde H_1$. We estimate $S_{B_3} (N)$ using the first part of Lemma~\ref{Cao-Zhai}:
\[
	\mathcal{B}_2 \ll_{\varepsilon_1} L^{2\varepsilon_1} \Big( L \tilde H_1 \tilde H_2 + X^{-1} (L \tilde H_1 \tilde H_2)^2 + L^{-2} \tilde H_1^2 \tilde H_2^6 + \tilde H_1^2 \tilde H_2^{8/3} \Big).
\] We can similarly conclude that the expression $X^{1/8} \mathcal{B}_1^{1/8} \mathcal{B}_2^{1/8}$ is maximized when $\tilde H_1$ and $\tilde H_2$ are both maximal, which occurs when $\tilde H_1 = H_1$ and $\tilde H_2 = L^{\varepsilon_1} H_1$. Then
\begin{multline*}
	S_{B_3} (N) \ll \frac{N^{\varepsilon + \varepsilon_1}}{N^2} \sum_{K, Y_1, Y_2} K^{(2-\theta) / (1-\theta) - 1} \max_{K \le \tilde K \le 2K} \sum_R K^{-3/(2(1-\theta))} (Y_1 Y_2)^{1-\theta/2} \cdot \\
	\frac{\bigl( \# \mathcal{N}_R (K, Y_1, Y_2) L \bigr)^{1/4}}{Z^{1/2} (H_1 H_2)^{1/4}} \bigl( Z^{1-1/\theta} L^{1/\theta - 2} H_1^2 \bigr)^{1/8} \cdot \\
	\Big( (KY_2^{\theta-1})^2 + \bigl( (KY_2^{\theta-1})^{\theta / (1-\theta)} Z^{1/\theta} L^{2 - 1\theta} H_1^{-2} \bigr) (KY_2^{\theta-1})^4 \Big)^{1/8} \cdot \\
	\Big( L H_1^2 +  \frac{(L H_1^2)^2}{Z^{1-1/\theta} L^{1/\theta-2} H_1^2} + L^{-2} H_1^8 + H_1^{14/3} \Big)^{1/8}.
\end{multline*}

\textit{Case 1. $Y_2 / 2 \le Y_1 \le 2Y_2$.} \\

Similarly to~\eqref{S_B_3_1_before_R}, we have
\begin{multline*} 
	S_{B_3}^{(1)} (N) \ll \frac{N^{\varepsilon + \varepsilon_1}}{N^2} \sum_{\substack{K, Y_1, Y_2 \\ Y_2 / 2 \le Y_1 \le 2Y_2}} Y_1^{1 - \theta / 2} Y_2^{1/2 + \theta / 4} \sum_R \frac{R}{(H_1 H_2)^{1/4}} \cdot \\
	\bigl( K^2 Y_2^{-1} R^{-1} H_1^2 \bigr)^{1/8} \Big( (KY_2^{\theta-1})^2 + \frac{(KY_2^{\theta-1})^4}{K^3 Y_2^{\theta - 2} R^{-2} H_1^2} \Big)^{1/8} \cdot \\
	\Big( RY_2 K^{-1} H_1^2 + \frac{\bigl( RY_2 K^{-1} H_1^2 \bigr)^2}{K^2 Y_2^{-1} R^{-1} H_1^2} + \frac{H_1^8}{\bigl( RY_2 K^{-1} \bigr)^2} + H_1^{14/3} \Big)^{1/8}.
\end{multline*}

Again, the most contribution to the sum over $R$ comes from $R$ maximal, so we get
\begin{multline*}
	S_{B_3}^{(1)} (N) \ll \frac{N^{\varepsilon + \varepsilon_1}}{N^2} \sum_{\substack{K, Y_1, Y_2 \\ Y_2 /2 \le Y_1 \le 2Y_2}} \frac{KY_2^{1/2 + 3\theta/4}}{(H_1 H_2)^{1/4}} \cdot \\
	\bigl( KY_2^{-\theta} H_1^2 \bigr)^{1/8} \Big( (KY_2^{\theta-1})^2 + \frac{(KY_2^{\theta-1})^4}{KY_2^{-\theta} H_1^2} \Big)^{1/8} \cdot \\
	\Big( Y_2^{\theta} H_1^2 + \frac{(Y_2^{\theta}H_1^2)^2}{KY_2^{-\theta} H_1^2} + \frac{H_1^8}{Y_2^{2\theta}} + H_1^{14/3} \Big)^{1/8}.
\end{multline*}	

Next, we take maximum values of $K, Y_1, Y_2$, yielding
\begin{multline} \label{S_B_3_1_new}
	S_{B_3}^{(1)} (N) \ll N^{3\varepsilon + \varepsilon_1} \frac{N^{-1/2 + 3\theta / 4}}{(H_1 H_2)^{1/4}} \bigl( N^{1-\theta} H_1^2 \bigr)^{1/8} \Big( N^{2\theta} + \frac{N^{4\theta}}{N^{1-\theta} H_1^2} \Big)^{1/8} \cdot \\
	\Big( N^{\theta} H_1^2 + \frac{\bigl( N^{\theta} H_1^2 \bigr)^2}{N^{1-\theta} H_1^2 } + \frac{H_1^8}{N^{2\theta}} + H_1^{14/3} \Big)^{1/8}.
\end{multline}

\textit{Case 2. $Y_1 > 2Y_2$.} \\

We have
\[
	S_{B_3}^{(2)} (N) \ll \frac{\log N}{N^2} \sum_{\substack{K, Y_1, Y_2 \\ Y_1 > 2Y_2}} \frac{K Y_1^{1/4 + \theta/2} Y_2^{1/4 + \theta / 4}}{(H_1 H_2)^{1/4}} \bigl( KY_1^{-\theta} H_1^2 \bigr)^{1/8} \mathcal{B}_1^{1/8} \mathcal{B}_2^{1/8}.
\] 

\textit{Case 2.1. $(KY_2^{\theta-1})^{2+\varepsilon_1} \gg \bigl( (KY_2^{\theta-1})^{\theta / (1-\theta)} ZX^{-1} \bigr) (KY_2^{\theta-1})^{4+\varepsilon_1}$.} \\

Note that now $X \asymp KY_1^{-\theta} H_1^2 L^{O(\varepsilon_1)}$. Similarly to~\eqref{S_B_3_2_1_before_KY}, we get
\begin{multline*}
	S_{B_3}^{(2.1)} (N) \ll
	\frac{N^{\varepsilon + \varepsilon_1}}{N^2} \sum_{\substack{K, Y_1, Y_2 \\ Y_1 > 2Y_2 \\ Y_2^{2-\theta} \gg KY_1^{2\theta} (H_1^2)^{-1} N^{O(\varepsilon_1)}}} \frac{K Y_1^{1/4 + \theta/2} Y_2^{1/4 + \theta / 4}}{(H_1 H_2)^{1/4}} \cdot \\
	\bigl( KY_1^{-\theta} H_1^2 \bigr)^{1/8} (KY_2^{\theta-1})^{1/4} 
	\Big( Y_1^{\theta} H_1^2 + \frac{(Y_1^{\theta} H_1^2)^2}{KY_1^{-\theta} H_1^2} + \frac{H_1^8}{Y_1^{2\theta}} + H_1^{14/3} \Big)^{1/8}.
\end{multline*} Taking the maximal $K, Y_1, Y_2$, we obtain
\begin{multline} \label{S_B_3_2_1_new}
	S_{B_3}^{(2.1)} (N) \ll N^{3\varepsilon + \varepsilon_1} \frac{N^{-1/2 + 3\theta / 4}}{(H_1 H_2)^{1/4}} \bigl( N^{1-\theta} H_1^2 \bigr)^{1/8} N^{\theta / 4} \cdot \\
	\Big( N^{\theta} H_1^2 + \frac{(N^{\theta} H_1^2)^2}{N^{1-\theta} H_1^2} + \frac{H_1^8}{N^{2\theta}} + H_1^{14/3} \Big)^{1/8}.
\end{multline} From~\eqref{S_B_3_1_new} and~\eqref{S_B_3_2_1_new}, we find
\begin{multline*}
	S_{B_3}^{(1)} (N) + S_{B_3}^{(2.1)} (N) \ll N^{3\varepsilon + \varepsilon_1} \frac{N^{-1/2 + 3\theta / 4}}{(H_1 H_2)^{1/4}} \bigl( N^{1-\theta} H_1^2 \bigr)^{1/8} \Big( N^{2\theta} + \frac{N^{4\theta}}{N^{1-\theta} H_1^2} \Big)^{1/8} \cdot \\
	\Big( N^{\theta} H_1^2 + \frac{\bigl( N^{\theta} H_1^2 \bigr)^2}{N^{1-\theta} H_1^2 } + \frac{H_1^8}{N^{2\theta}} + H_1^{14/3} \Big)^{1/8}.
\end{multline*} With our choice of $(\kappa_1, \kappa_2)$ we have $N^{2\theta} > N^{4\theta} / (N^{1-\theta} H_1^2)$. Thus,
\begin{multline*}
	S_{B_3}^{(1)} (N) + S_{B_3}^{(2.1)} (N) \ll N^{3\varepsilon + \varepsilon_1} \frac{N^{-3/8 + 7\theta / 8}}{H_2^{1/4}} \cdot \\
	\Big( N^{\theta} H_1^2 + \frac{\bigl( N^{\theta} H_1^2 \bigr)^2}{N^{1-\theta} H_1^2 } + \frac{H_1^8}{N^{2\theta}} + H_1^{14/3} \Big)^{1/8} = N^{4\varepsilon} \Big(
	N^{-3/8 + \theta + \kappa_1 \theta / 4 - \kappa_2 \theta / 4} + \\
	N^{-1/2 + 5\theta / 4 + \kappa_1 \theta / 4 - \kappa_2 \theta / 4} +
	N^{-3/8 + 5\theta / 8 + \kappa_1 \theta - \kappa_2 \theta / 4} +
	N^{-3/8 + 7\theta / 8 + 7\kappa_1 \theta / 12 - \kappa_2 \theta / 4} \Big),
\end{multline*} which is $o(1)$ when
\[
	\theta < \min \Big( \frac{3}{8 + 2\kappa_1 - 2\kappa_2}, \
	\frac{2}{5 + \kappa_1 - \kappa_2}, \
	\frac{3}{5 + 8\kappa_1 - 2\kappa_2}, \
	\frac{9}{21 + 14\kappa_1 - 6\kappa_2} \Big) - 5\varepsilon,
\] which covers the range $[0, \frac{43}{117})$ when $(\kappa_1, \kappa_2) = (\frac{12}{43}, \frac{24}{43})$. \\

\textit{Case 2.2. $(KY_2^{\theta-1})^{2+\varepsilon_1} \ll \bigl( (KY_2^{\theta-1})^{\theta / (1-\theta)} ZX^{-1} \bigr) (KY_2^{\theta-1})^{4+\varepsilon_1}$.} \\

Similarly to~\eqref{S_B_3_2_2_first}, we find
\begin{multline*}
	S_{B_3}^{(2.2)} (N) \ll \frac{N^{\varepsilon + \varepsilon_1}}{N^2} \sum_{\substack{K, Y_1, Y_2 \\ Y_1 > 2Y_2 \\ Y_2^{2-\theta} \ll KY_1^{2\theta} (H_1^2)^{-2} N^{O(\varepsilon_1)}}} \frac{KY_1^{1/4 + \theta / 2} Y_2^{1/4 + \theta / 4}}{(H_1 H_2)^{1/4}} \bigl( KY_1^{-\theta} H_1^2 \bigr)^{1/8} \cdot \\
	\Big( \frac{(KY_2^{\theta-1})^4}{KY_1^{-2\theta} Y_2^{\theta} H_1^2}  \Big)^{1/8} 
	\Big( Y_1^{\theta} H_1^2 + \frac{(Y_1^{\theta} H_1^2)^2}{KY_1^{-\theta} H_1^2} + \frac{H_1^8}{Y_1^{2\theta}} + H_1^{14/3} \Big)^{1/8}.
\end{multline*} The maximum contribution comes from $K \asymp N^{1+\varepsilon}$, $Y_1 \asymp N$, $Y_2 \asymp 1$:
\begin{multline*}
	S_{B_3}^{(2.2)} (N) \ll N^{3\varepsilon + \varepsilon_1} \frac{N^{-1/4 + 5\theta / 8}}{(H_1 H_2)^{1/4}} \Big( N^{\theta} H_1^2 + N^{3\theta-1} H_1^2 + \frac{H_1^8}{N^{2\theta}} + H_1^{14/3} \Big)^{1/8} \ll \\
	N^{4\varepsilon} \Big( 
	N^{-1/4 + 3\theta/4 - \kappa_2 \theta / 4} +
	N^{-3/8 + \theta - \kappa_2 \theta / 4} + \\
	N^{-1/4 + 3\theta / 8 + 3\kappa_1 \theta / 4 - \kappa_2 \theta / 4} +
	N^{-1/4 + 5\theta / 8 + \kappa_1 \theta / 3 - \kappa_2 \theta / 4} \Big),
\end{multline*} which is $o(1)$ when
\[
	\theta < \min \Big( \frac{1}{3 - \kappa_2}, \
	\frac{3}{8 - 2\kappa_2}, \
	\frac{2}{3 + 6\kappa_1 - 2\kappa_2}, \
	\frac{6}{15 + 8\kappa_1 - 6\kappa_2} \Big) - 5\varepsilon,
\] which covers the range $[0, \frac{43}{117})$ when $(\kappa_1, \kappa_2) = (\frac{12}{43}, \frac{24}{43})$.

\bibliographystyle{plain}
\bibliography{Bib}

\end{document}